\documentclass[12pt,fullpage]{article}

\usepackage{amsfonts}
\usepackage{xr}
\usepackage{graphicx}
\usepackage{amsmath}
\usepackage{epsfig}

\def\R{\mathbb{R}}
\def\N{\mathbb{N}}
\def\epsilon{\varepsilon}
\def\trait (#1) (#2) (#3){\vrule width #1pt height #2pt depth #3pt}
\def\fin{\hfill\trait (0.1) (5) (0) \trait (5) (0.1) (0) \kern-5pt 
\trait (5) (5) (-4.9) \trait (0.1) (5) (0)}

\newcommand{\SE}{\setcounter{equation}{0} \section}
\newcommand{\be}{\begin{equation}}
\newcommand{\ee}{\end{equation}}
\newcommand{\baa}{\begin{array}}
\newcommand{\eaa}{\end{array}}
\newcommand{\ba}{\begin{eqnarray}}
\newcommand{\ea}{\end{eqnarray}}

\newtheorem{theo}{\bf Theorem}[section]
\newtheorem{lem}[theo]{\bf Lemma}
\newtheorem{pro}[theo]{\bf Proposition}
\newtheorem{cor}[theo]{\bf Corollary}

\newtheorem{rem}[theo]{\bf Remark}

\begin{document}
\date{}
\title{\bf{A Faber-Krahn inequality with drift}}
\author{Fran\c cois Hamel$^{\hbox{\small{ a}}}$, Nikolai 
Nadirashvili$^{\hbox{\small{ b}}}$ and Emmanuel Russ$^{\hbox{\small{ 
a}}}$\\
\\
\footnotesize{$^{\hbox{a }}$Universit\'e Aix-Marseille III, LATP, 
Facult\'e des Sciences et Techniques, Case cour A}\\
\footnotesize{Avenue Escadrille Normandie-Niemen, F-13397 Marseille 
Cedex 20, France}\\
\footnotesize{$^{\hbox{b }}$ CNRS, LATP, CMI, 39 rue F. Joliot-Curie, 
F-13453 Marseille Cedex 13, France}\\
\footnotesize{francois.hamel@univ.u-3mrs.fr, nicolas@cmi.univ-mrs.fr, 
emmanuel.russ@univ.u-3mrs.fr}}

\maketitle

\begin{abstract} Let $\Omega$ be a bounded $C^{2,\alpha}$ domain in $\R^n$ ($n\geq 1$, $0<\alpha<1$), $\Omega^{\ast}$ be the open Euclidean ball centered at $0$ having the same Lebesgue measure as $\Omega$, $\tau\geq 0$ and $v\in L^{\infty}(\Omega,\R^n)$ with $\left\Vert v\right\Vert_{\infty}\leq \tau$. If $\lambda_{1}(\Omega,\tau)$ denotes the principal eigenvalue of the operator $-\Delta+v\cdot\nabla$ in $\Omega$ with Dirichlet boundary condition, we establish that $\lambda_{1}(\Omega,v)\geq \lambda_{1}(\Omega^{\ast},\tau e_{r})$ where $e_{r}(x)=x/\left\vert x\right\vert$. Moreover, equality holds only when, up to translation, $\Omega=\Omega^{\ast}$ and $v=\tau e_{r}$. This result can be viewed as an isoperimetric inequality for the first eigenvalue of the Dirichlet Laplacian with drift. It generalizes the celebrated Rayleigh-Faber-Krahn inequality for the first eigenvalue of the Dirichlet Laplacian.
\end{abstract}


\SE{Introduction and main results}\label{intro}

Throughout all the paper, $n\geq 1$ denotes an integer in $\N^*=\N\backslash\{0\}$. By ``domain'', we mean an open connected subset of $\R^n$, and we denote by ${\mathcal C}$ the set of all bounded domains of $\R^n$ which are of class $C^{2,\alpha}$ for some $0<\alpha<1$. For any measurable subset $A\subset\R^n$, $\left\vert A\right\vert$ stands for the standard $n$-dimensional Lebesgue measure of $A$. Throughout the paper, $B^n_r$ denotes the open Euclidean ball of $\R^n$ with center $0$ and radius $r>0$, and we set $\alpha_n=|B^n_1|=\pi^{n/2}/\Gamma(n/2+1)$. For $\Omega\in{\mathcal C}$, we define $\Omega^{\ast}$ as the ball $B^n_{(|\Omega|/\alpha_n)^{1/n}}$ having the same measure as $\Omega$. Finally, if $\Omega\in {\mathcal C}$ and $v:\Omega\rightarrow \R^n$ is measurable, $\left\vert v\right\vert$ will denote the Euclidean norm of $v$, and we say that $v\in L^{\infty}(\Omega,\R^n)$ if $\left\vert v\right\vert\in L^{\infty}$, and write (somewhat abusively) $\left\Vert v\right\Vert_{\infty}$ instead of $\left\Vert\ \left\vert v\right\vert\ \right\Vert_{L^{\infty}(\Omega)}$. 

If $\lambda_1(\Omega)$ denotes the first eigenvalue of the Laplace operator $-\Delta$ with Dirichlet boundary condition (Dirichlet Laplacian), in an open bounded smooth set $\Omega\subset\R^n$, it is well-known that $\lambda_1(\Omega)\geq \lambda_1(\Omega^{\ast})$ and that the inequality is strict unless $\Omega$ is a ball. Since $\lambda_1(\Omega^{\ast})$ can be explicitly computed, this result provides the classical Rayleigh-Faber-Krahn inequality, which states that
\begin{equation} \label{RFK}
\lambda_1(\Omega)\geq\vert\Omega\vert^{-2/n}\alpha_n^{2/n}j_{n/2-1,1}^2,
\end{equation}
where $j_{m,1}$ the first positive zero of the Bessel function $J_m$. Moreover, equality in (\ref{RFK}) is attained if and only if $\Omega$ is a ball. This result was first conjectured by Rayleigh (1894/1896) for $n=2$ (\cite{Rayleigh} vol. I, pp. 339-345), and proved independently by Faber (\cite{Faber}, 1923) and Krahn (\cite{Krahndim2}, 1925) for $n=2$, and by Krahn for all $n$ in \cite{Krahndimn} (1926; see \cite{Krahndimnenglish} for the English translation). 

Many other optimization results for the eigenvalues of the Dirichlet Laplacian have been proved. For instance, the minimum of $\lambda_2(\Omega)$ among open bounded sets $\Omega\subset\R^n$ with given Lebesgue measure is achieved by the union of two identical balls (this result is attributed to Szeg\"o, see \cite{polya}). Very few things seem to be known about optimization problems for the other eigenvalues, see \cite{bucurhenrot,henrot,polya,polya2,wk}.

Various optimization results are known about functions of the eigenvalues. For instance, the Payne-P\'olya-Weinberger conjecture (see \cite{ppw}) on the ratio of 
the first two eigenvalues was proved by Ashbaugh and Benguria (\cite{ashben}), in any dimension $n$: namely, for any bounded domain $\Omega\subset \R^n$, 
$\lambda_2(\Omega)/\lambda_1(\Omega)\leq\lambda_2(\Omega^{\ast})/\lambda_1(\Omega^{\ast})$, and the equality is attained only when $\Omega$ is a ball. 
The same result was also extended in \cite{ashben} for elliptic operators in divergence form with definite weight. 
We also refer to \cite{ab2,ab3,ahs,barnes,co,karaa,karaa2,ly,m,ppw,polya2} for further bounds or other optimization results for some eigenvalues or some functions of the eigenvalues in fixed or varying domains.

Other boundary conditions may also be considered. For instance, if $\mu_2(\Omega)$ is the first non-trivial eigenvalue of the Laplacian under Neumann boundary condition, $\mu_2(\Omega)\leq\mu_2(\Omega^{\ast})$ and the equality is attained only when $\Omega$ is a ball (see \cite{sz} in dimension $n=2$, and \cite{wein} in any dimension). Bounds or optimization results for other eigenvalues of the Laplacian under Neumann boundary condition (\cite{ppw,polya2,sz,wein}, see also \cite{bandle1} for inhomogeneous problems), for Robin boundary condition (\cite{bossel}) or for the Stekloff eigenvalue problem (\cite{brock}) have also been established.

We also mention another Rayleigh conjecture for the lowest eigenvalue of the clamped plate. If $\Omega$ is a smooth bounded open subset of $\R^2$, denote by $\Lambda_1(\Omega)$ the lowest eigenvalue of the operator $\Delta^2$, so that $\Delta^2u_1=\Lambda_1(\Omega)u_1$ in $\Omega$ with $u_1=|\nabla u_1|=0$ on $\partial\Omega$ and $u_1>0$ in $\Omega$. The second author proved in \cite{nad} that $\Lambda_1(\Omega)\geq \Lambda_1(\Omega^{\ast})$ and that equality holds only when $\Omega$ is a ball (disk, in dimension 2). The analogous result was also established in $\R^3$ in \cite{ashbenclamped}, whereas the problem is still open in higher dimensions. 

Very nice and much more complete surveys of all these topics and additional results can be found in \cite{bandle2,henrot,o} and the references therein.

All above problems concern self-adjoint operators. In the present paper, we focus on optimization problems for the first eigenvalue of the {\it non-self-adjoint} Laplace operator {\it with a drift} under Dirichlet boundary condition.

For any domain $\Omega\in{\mathcal{C}}$, for any $v\in L^{\infty}(\Omega,\R^n)$, we call $\lambda_1(\Omega,v)$ the first eigenvalue of $L=-\Delta+v\cdot\nabla$ with Dirichlet boundary condition on $\partial\Omega$, and $\varphi_{\Omega,v}$ the corresponding (unique) positive eigenfunction with $L^{\infty}$-norm equal to $1$. In the sequel, $\varphi_{\Omega,v}$ will be denoted by $\varphi$ when the context makes clear what $\Omega$ and $v$ are. Recall that the maximum principle holds for $L$, and, as a consequence, $\lambda_1(\Omega,v)>0$ (see \cite{bnv}). One has
\begin{equation} \label{eq}
\left\{
\begin{array}{l}
-\Delta\varphi+v\cdot\nabla\varphi=\lambda_1(\Omega,v)\varphi\mbox{ 
in }\Omega,\\
\\
\varphi>0\mbox{ in }\Omega,\ \varphi=0\mbox{ on }\partial\Omega,\ 
\left\Vert \varphi\right\Vert_{L^{\infty}(\Omega)}=1.
\end{array}
\right.
\end{equation}
Moreover, by standard elliptic estimates (see \cite{adn,gt}), $\varphi\in W^{2,p}(\Omega)$ for all $1\leq p<+\infty$, whence, up to the choice of the continuous representant in the class of $\varphi$, $\varphi\in C^{1,\beta}(\overline{\Omega})$ for all $0\leq\beta<1$. Recall also that if $\lambda$ is any eigenvalue for the operator $L$, then either $\lambda=\lambda_1(\Omega,v)$ or $\mbox{Re}(\lambda)>\lambda_1(\Omega,v)$, and that, if $\psi$ is any positive eigenfunction in $\Omega$ for $L$ corresponding to the eigenvalue $\lambda$, then actually $\lambda=\lambda_1(\Omega,v)$ and $\psi$ and $\varphi$ are proportional (see \cite{bnv} again).

For all $x\neq 0$, set
\[
e_r(x)=\frac{x}{\left\vert x\right\vert}. 
\]

Our main result is the following one:

\begin{theo} \label{th1}
For any dimension $n\ge 1$, any $\Omega\in {\mathcal C}$, any $\tau\geq 0$ and any $v\in L^{\infty}(\Omega,\R^n)$ satisfying $\left\Vert v\right\Vert_{\infty}\leq \tau$, 
\begin{equation} \label{main}
\lambda_1(\Omega,v)\geq \lambda_1(\Omega^{\ast},\tau e_r).
\end{equation}
Moreover, equality holds only when, up to translation, $\Omega=\Omega^{\ast}$ and $v=\tau e_r$, namely when there exists $x_0\in\R^n$ such that $(\Omega,v)=(x_0+\Omega^*,\tau e_r(\cdot-x_0))$.
\end{theo} 

Theorem \ref{th1} can then be viewed as a natural extension of the first Rayleigh conjecture to the Dirichlet Laplacian with drift in any dimension $n$.

A rough parabolic interpretation of Theorem \ref{th1} can be the following one: consider the evolution equation $u_t=\Delta u-v\cdot\nabla u$ in $\Omega$, 
for $t>0$, with Dirichlet boundary condition on $\partial\Omega$, and with an initial datum at $t=0$. 
Roughly speaking, minimizing $\lambda_1(\Omega,v)$ (with given $|\Omega|$ and with $\|v\|_{\infty}\le\tau$ can be interpreted as looking for the slowest 
exponential time-decay of the solution $u$. The best way to do that is to try to minimize the boundary effects, namely to have the domain as round as possible, and it is not unreasonable to say that the vector field $-v$ should as much as possible point inwards the domain to avoid the drift towards the boundary. Of course diffusion, boundary losses and transport phenomena take place simultaneously, but these heuristic arguments tend to lead to the optimal couple $(\Omega,-v)=(\Omega^*,-\tau e_r)$ (up to translation).

As a corollary of Theorem \ref{th1}, we obtain the following Faber-Krahn type inequality for the Dirichlet Laplacian with drift, which extends the classical Rayleigh-Faber-Krahn inequality for the Dirichlet Laplacian. Namely we get a lower bound for $\lambda_1(\Omega,v)$, which depends only on $\left\vert \Omega\right\vert$ and $\left\Vert v\right\Vert_{\infty}$:

\begin{cor}
For each $n\ge 1$, there exists a function $F_n:(0,+\infty)\times [0,+\infty)\rightarrow (0,+\infty)$, defined by $F_n(m,\tau)=\lambda_1(B^n_{(m/\alpha_n)^{1/n}},\tau e_r)$ such that, for any domain $\Omega\in {\mathcal C}$ and any $v\in L^{\infty}(\Omega,\R^n)$,
\be\label{fk}
\lambda_1(\Omega,v)\geq F_n(|\Omega|,\|v\|_{\infty})
\ee
and equality holds if and only if, up to translation, $\Omega=\Omega^*$ and $v=\|v\|_{\infty}e_r$.
\end{cor}

\begin{rem} \label{rem1} {\rm For fixed $n\in\N^*$, $m>0$ and $\tau\geq 0$, inequality (\ref{main}) of Theorem \ref{th1} may be reformulated in the following way: $\displaystyle{\mathop{\inf}_{\Omega\in{\mathcal{C}},\ \left\vert \Omega\right\vert=m,\ \left\Vert v\right\Vert_{\infty}\leq \tau}}\lambda_1(\Omega,v)=\lambda_1(B,\tau e_r)$, where $B=B^n_{(m/\alpha_n)^{1/n}}$. Since $\lambda_1\left(\Omega',\left. v\right\vert_{\Omega'}\right)>\lambda_1(\Omega,v)$ for all $\Omega$, $\Omega'\in{\mathcal{C}}$, $\displaystyle \Omega'\underset{\neq}\subset\Omega$ and $v\in L^{\infty}(\Omega,\R^n)$ (see \cite{bnv}), Theorem \ref{th1} yields at once
$$\inf_{\Omega\in{\mathcal{C}},\ \left\vert \Omega\right\vert\leq m,\ \left\Vert v\right\Vert_{\infty}\leq \tau}\lambda_1(\Omega,v)=\lambda_1(B,\tau e_r),$$
and the infimum is reached for and only for $\Omega=B$ and $v=\tau e_r$ (up to translation).\par
The above inequalities (\ref{main}) and (\ref{fk}) can also be generalized to more general open sets $\Omega$. First, one has $\displaystyle{\mathop{\inf}_{\Omega\in{\mathcal{C}}',\ \left\vert \Omega\right\vert\leq m,\ \left\Vert v\right\Vert_{\infty}\leq \tau}}\lambda_1(\Omega,v)=\lambda_1(B,\tau e_r)$, where ${\mathcal{C}}'$ denotes the set of all open bounded subsets $\Omega$ of $\R^n$ of class $C^{2,\alpha}$ for some $0<\alpha<1$, with finite number (maybe not reduced to $1$) of connected components (for $\Omega\in{\mathcal{C}}'$ and $v\in L^{\infty}(\Omega,\R^n)$, $\lambda_1(\Omega,v)$ is the minimum, over $k$, of the eigenvalues $\lambda_1(\Omega_k,v\vert_{\Omega_k})$, where the $\Omega_k$'s are the connected components of $\Omega$). Thus, inequalities (\ref{main}) and (\ref{fk}) hold for $\Omega\in{\mathcal{C}}'$ and the case of equality holds only if, up to translation, $\Omega=B$ and $v=\tau e_r$.\par
Furthermore, for non-smooth and possibly unbounded $\Omega$ with finite measure, following \cite{bnv}, one can still define $\lambda_1(\Omega,v)$, 
as $\lambda_1(\Omega,v)=\displaystyle{\mathop{\inf}_{\Omega'\subset\Omega,\ \Omega'\in{\mathcal{C'}}}}\lambda_1(\Omega',v|_{\Omega'})$. Since $F_n(m,\tau)$ is 
decreasing in both $m>0$ and $\tau\ge 0$ (see Remark \ref{contm}), inequalities 
(\ref{main}) and (\ref{fk}) still hold.}
\end{rem}

\begin{rem} \label{rem2}{\rm Let us now discuss the behavior of $F_n(m,\tau)$ for large $\tau$ (see Section \ref{secFn} for details). First, for all $m>0$, $\tau^{-2}e^{\tau m/2}F_1(m,\tau)\to 1$ as $\tau+\infty$, and one even has
\be\label{F1}
\exists\ C(m)\ge 0,\ \exists\ \tau_0\ge 0,\ \forall\ \tau\ge\tau_0,\quad|\tau^{-2}e^{\tau m/2}F_1(m,\tau)-1|\le C(m)\tau e^{-\tau m/2}.
\ee
Moreover, for all $n\geq 2$ and $m>0$, $F_n(m,\tau)>F_1(2(m/\alpha_n)^{1/n},\tau)$ for all $\tau\ge 0$, and
\be\label{Fntau}
-\tau^{-1}\log F_n(m,\tau)\to m^{1/n}\alpha_n^{-1/n}\hbox{ as }\tau\rightarrow +\infty.
\ee
In \cite{fr}, with probabilistic arguments, Friedman proved some lower and upper logarithmic estimates, as $\varepsilon\to 0^+$, for the first eigenvalue of general elliptic operators $-a_{ij}\varepsilon^2\partial_{ij}+b_i\partial_i$ with $C^1$ drifts $-b=-(b_1,\ldots,b_n)$ pointing inwards on the boundary. Apart from the fact that the vector field $e_r$ is not $C^1$ at the origin, the general result of Friedman would imply the asymptotics (\ref{Fntau}) for $\log F_n(m,\tau)=\log \lambda_1(B^n_{(m/\alpha_n)^{1/n}},\tau e_r)$. For the sake of completeness, we give in Appendix (Section \ref{secFn}) a proof of (\ref{Fntau}) with elementary analytic arguments. There, we also prove the precise equivalent of $F_1(m,\tau)$ for large $\tau$. However, giving an equivalent for $F_n(m,\tau)$ when $\tau$ is large and $n\geq 2$ is an open question.}
\end{rem} 

The first step in the proof of Theorem \ref{th1}, which has its own interest, is the optimization of $\lambda_1(\Omega,v)$ when $\Omega$ is a fixed domain and the $L^{\infty}$ norm of $v$ is controlled. Namely, for any $\tau\geq 0$, set
\[
\underline{\lambda}(\Omega,\tau)=\inf_{\left\Vert v\right\Vert_{\infty}\leq\tau} \lambda_1(\Omega,v)\mbox{ and } \overline{\lambda}(\Omega,\tau)=\sup_{\left\Vert v\right\Vert_{\infty}\leq\tau} \lambda_1(\Omega,v).
\]
It turns out that this optimization problem also has a unique solution:

\begin{theo} \label{th2}
Let $\Omega$ be a domain in ${\mathcal C}$ $($of class $C^{2,\alpha}$ for some $0<\alpha<1)$ and let $\tau\geq 0$ be fixed.
\begin{itemize}
\item[$(a)$] 
There exists a unique vector field $\underline{v}\in L^{\infty}(\Omega)$ with $\left\Vert \underline{v}\right\Vert_{\infty}\leq \tau$ such that $\underline{\lambda}(\Omega,\tau)=\lambda_1(\Omega,\underline{v})$ $(>0)$, and this field satisfies $|\underline{v}(x)|=\tau$ almost everywhere in $\Omega$. Moreover, the corresponding principal eigenfunction $\underline{\varphi}=\varphi_{\Omega,\underline{v}}$ is of class $C^{2,\alpha}(\overline{\Omega})$ and $\underline{v}\cdot \nabla\underline{\varphi}=-\tau\left\vert \nabla \underline{\varphi}\right\vert$ almost everywhere in $\Omega$. The function $\underline{\varphi}$ is then a solution of the following nonlinear problem
\begin{equation} \label{eqmin}
-\Delta\underline{\varphi}-\tau\left\vert \nabla \underline{\varphi}\right\vert=\underline{\lambda}(\Omega,\tau)\underline{\varphi}\mbox{ in }\Omega.
\end{equation}
Moreover, if $\psi$ is a function of class $C^{2,\alpha}(\overline{\Omega})$ such that $\psi>0$ in $\Omega$, $\psi=0$ on $\partial\Omega$, $\left\Vert \psi\right\Vert_{\infty}=1$ and if $\mu\in\R$ is such that $(\ref{eqmin})$ holds with $\psi$ and $\mu$ instead of $\underline{\varphi}$ and $\underline{\lambda}(\Omega,\tau)$, then $\psi=\underline{\varphi}$ and $\mu=\underline{\lambda}(\Omega,\tau)$. 
\item[$(b)$]
There exists a unique vector field $\overline{v}\in L^{\infty}(\Omega)$ with $\left\Vert \overline{v}\right\Vert_{\infty}\leq \tau$ such that $\overline{\lambda}(\Omega,\tau)=\lambda_1(\Omega,\overline{v})$ $(>0)$, and this field satisfies $|\overline{v}(x)|=\tau$ almost everywhere in $\Omega$. Moreover, the corresponding principal eigenfunction $\overline{\varphi}=\varphi_{\Omega,\overline{v}}$ is of class $C^{2,\alpha}(\overline{\Omega})$ and $\overline{v}\cdot \nabla\overline{\varphi}=\tau\left\vert \nabla \overline{\varphi}\right\vert$ almost everywhere in $\Omega$. The function $\overline{\varphi}$ is then a solution of the following nonlinear problem
\begin{equation} \label{eqmax}
-\Delta\overline{\varphi}+\tau\left\vert \nabla \overline{\varphi}\right\vert=\overline{\lambda}(\Omega,\tau)\overline{\varphi}\mbox{ in }\Omega.
\end{equation}
Moreover, if $\psi$ is a function of class $C^{2,\alpha}(\overline{\Omega})$ such that $\psi>0$ in $\Omega$, $\psi=0$ on $\partial\Omega$, $\left\Vert \psi\right\Vert_{\infty}=1$ and if $\mu\in\R$ is such that $(\ref{eqmax})$ holds with $\psi$ and $\mu$ instead of $\overline{\varphi}$ and $\overline{\lambda}(\Omega,\tau)$, then $\psi=\overline{\varphi}$ and $\mu=\overline{\lambda}(\Omega,\tau)$.
\end{itemize}
\end{theo}

When $\Omega$ is a ball (up to translation, one assumes that it is centered at the origin), one can provide an explicit expression of $\underline{v}$ and $\overline{v}$:

\begin{theo} \label{th3}
Assume that $\Omega=B=B^n_R$ for some radius $R>0$, and let $\tau\ge 0$ be fixed. Then $\underline{v}=\tau e_r$ and $\overline{v}=-\tau e_r$ where $\underline{v}$ and $\overline{v}$ are defined in Theorem \ref{th2}. One therefore has:
\begin{equation} \label{eqball}\left\{\begin{array}{rcll}
-\Delta\underline{\varphi}+\tau e_{r}\cdot\nabla\underline{\varphi} & = & \underline{\lambda}(\Omega,\tau)\underline{\varphi} &\mbox{ in }B,\\
\\
-\Delta\overline{\varphi}-\tau e_{r}\cdot\nabla\overline{\varphi} & = & \overline{\lambda}(\Omega,\tau)\overline{\varphi} &\mbox{ in }B.\end{array}\right.
\end{equation}
Moreover, the functions $\underline{\varphi}$ and $\overline{\varphi}$ are radially decreasing in $\overline{B}$, which means that there are two decreasing functions $\underline{\phi}$, $\overline{\phi}:[0,R]\to[0,+\infty)$ such that $\underline{\varphi}(x)=\underline{\phi}(|x|)$ and $\overline{\varphi}(x)=\overline{\phi}(|x|)$ for all $x\in\overline{B}$.
\end{theo}

\begin{rem}{\rm Notice that a corollary of Theorems \ref{th1}, \ref{th2} and \ref{th3} is that, for all $\tau\ge 0$ and $\Omega\in{\mathcal{C}}$, $\underline{\lambda}(\Omega,\tau)\ge\underline{\lambda}(\Omega^*,\tau)$, and equality holds only when $\Omega$ is a ball.}
\end{rem}

Let us now give a few additional comments about Theorems \ref{th1}, \ref{th2} and \ref{th3}. First, what happens for other optimization problems with analogous constraints? For a fixed domain $\Omega\in {\mathcal C}$, if we drop the condition $\left\Vert v\right\Vert_{\infty}\leq \tau$, one can prove that
\be\label{limits}
\inf_{v\in L^{\infty}(\Omega,\R^n)} \lambda_1(\Omega,v)=\inf_{\tau\ge 0}\ \underline{\lambda}(\Omega,\tau)=0,\ 
\sup_{v\in L^{\infty}(\Omega,\R^n)} \lambda_1(\Omega,v)=\sup_{\tau\ge 0}\ \overline{\lambda}(\Omega,\tau)=+\infty.
\ee
Actually, we prove in Lemmata \ref{continuity} and \ref{continuitybis} that the function $\tau\mapsto\underline{\lambda}(\Omega,\tau)$ (resp. $\tau\mapsto\overline{\lambda}(\Omega,\tau)$) is decreasing (resp. increasing) in $[0,+\infty)$. Therefore, (\ref{limits}) means that $\underline{\lambda}(\Omega,\tau)\to 0$ and $\overline{\lambda}(\Omega,\tau)\to+\infty$ as $\tau\to+\infty$. The proof of first assertion (about the infimum) follows at once from formula (\ref{Fntau}): indeed, choose a ball $B$ included in $\Omega$, call $x_0$ its center, and let $v\in L^{\infty}(\Omega,\R^n)$ be such that $v(x)=\tau e_r(x-x_0)$ for all $x\in B$; one then has $0<\lambda_1(\Omega,v)\le\lambda_1(B,v|_{B})=F_n(|B|,\tau)\to 0$ as $\tau\to+\infty$. For the proof of the other assertion, let $\tau>0$ and define $v=\tau e_1$, where, for all $x\in \R^n$, $e_1(x)=(1,0,\ldots,0)$. Straightforward computations show that the function $\psi(x)=e^{-\tau x_1/2}\varphi_{\Omega,\tau e_1}(x)$, which is positive in $\Omega$, satisfies $-\Delta \psi+(\tau^2/4)\psi=\lambda_1(\Omega,\tau e_1)\psi$ in $\Omega$. From the characterization of the first eigenvalue, one concludes that $\lambda_1(\Omega,\tau e_1)=\tau^2/4+\lambda_1(\Omega)$, where $\lambda_1(\Omega)$ is the first eigenvalue of the Dirichlet Laplacian in $\Omega$. This obviously implies the desired result.

Notice that, when $v\in L^{\infty}(\Omega)$ is divergence free (in the sense of distributions), $\lambda_1(\Omega,v)\ge\lambda_1(\Omega,0)=\lambda_1(\Omega)$ (indeed, multiply (\ref{eq}) by $\varphi$ and integrate by parts). Thus, it immediately follows that $\displaystyle{\mathop{\inf}_{\|v\|_{\infty}\le\tau,\ \hbox{\small{div}}(v)=0}}\lambda_1(\Omega,v)=\lambda_1(\Omega)$ for all $\Omega\in{\mathcal{C}}$ and $\tau\ge 0$. We also refer to \cite{bhn} for a detailed analysis of the behavior of $\lambda_1(\Omega,A\ v)$ when $A\to+\infty$ and $v$ is a fixed divergence free vector field in $L^{\infty}(\Omega)$.

Let now $v\in L^{\infty}(\R^n,\R^n)$ be fixed, call $\|v\|_{\infty}=\|\ |v|\ \|_{L^{\infty}(\R^n,\R^n)}$, let $\Omega$ vary and define
$$\underline{\underline{\lambda}}(v,m)=\inf_{\Omega\in{\mathcal{C}},\ |\Omega|=m}\lambda_1(\Omega,v|_{\Omega})=\inf_{\Omega\in{\mathcal{C}},\ |\Omega|\le m}\lambda_1(\Omega,v|_{\Omega})\ \hbox{ and }\ \overline{\overline{\lambda}}(v,m)=\sup_{\Omega\in{\mathcal{C}},\ |\Omega|=m}\lambda_1(\Omega,v|_{\Omega})$$
for each $m>0$. Because of (\ref{supomega}) below, there holds $\overline{\overline{\lambda}}(v,m)=+\infty$. On the other hand, $\underline{\underline{\lambda}}(v,m)\ge F_n(m,\|v\|_{\infty})$. However, unlike $\underline{\lambda}(\Omega,\tau)$ or $\overline{\lambda}(\Omega,\tau)$, the infimum in $\underline{\underline{\lambda}}(v,m)$ may not be reached: indeed, let $(x_k)_{k\ge k_0}$ be a sequence of points in $\R^n$ such that the balls $B_k=x_k+B^n_{(m/\alpha_n)^{1/n}}$ are pairwise disjoint and choose $k_0\in\N$ and $v\in L^{\infty}(\R^n)$ so that $2^{-k_0}\le\|v\|_{\infty}$, $v|_{B_k}=(\|v\|_{\infty}-2^{-k})\ e_r(\cdot-x_k)$ and $v=0$ outside the balls $B_k$; from Theorems \ref{th1}, \ref{th3} and Lemma \ref{continuity} below, one can easily check that, for each $\Omega\in{\mathcal{C}}$ with $|\Omega|\le m$, $\lambda_1(\Omega,v|_{\Omega})>\underline{\underline{\lambda}}(v,m)=F_n(m,\|v\|_{\infty})$.

Consider now optimization problems (different from the one solved in Theorem \ref{th1}) where both $\Omega$ and $v$ vary. Let $m>0$ and $\tau\geq 0$ be fixed. As was already mentioned in Remark~\ref{rem1},
\[
\inf_{\Omega\in{\mathcal{C}},\ \left\vert \Omega\right\vert\leq m,\ \left\Vert v\right\Vert_{\infty}\leq \tau}\lambda_1(\Omega,v)=\inf_{\Omega\in{\mathcal{C}},\ \left\vert \Omega\right\vert\leq m}\underline{\lambda}(\Omega,\tau)=\lambda_1(\Omega^{\ast},\tau e_r).
\]
The optimization problem for $\underline{\lambda}(\Omega,\tau)$ when $\left\vert \Omega\right\vert=m$ with a supremum instead of the infimum has the following solution: 

\begin{pro}\label{prosup} One has
\be\label{supomega}
\sup_{\Omega\in{\mathcal{C}},\ \left\vert \Omega\right\vert=m}\underline{\lambda}(\Omega,\tau)=+\infty,
\ee
whence $\displaystyle{\mathop{\sup}_{\Omega\in{\mathcal{C}},\ \left\vert\Omega\right\vert=m}}\overline{\lambda}(\Omega,\tau)=+\infty$.
\end{pro}

The proof follows from a min-max formula for $\lambda_1(\Omega,v)$ given in \cite{bnv} (see Section \ref{secminmax} in Appendix).\hfill\break

\noindent{\bf{Open problem.}} Finally, we mention as an open problem the characterization of
\[
\inf_{\Omega\in{\mathcal{C}},\ \left\vert \Omega\right\vert=m} \overline{\lambda}(\Omega,\tau)
\]
for given $m>0$ and $\tau\ge 0$. As for the infimum of $\underline{\lambda}(\Omega,\tau)$ with the constraint $|\Omega|=m$, is the infimum of $\overline{\lambda}(\Omega,\tau)$ achieved for the balls~?\hfill\break\par
 
We now turn to the strategy of the proof of our results. Remember first that the proof of the classical Rayleigh-Faber-Krahn inequality (\ref{RFK}) relies on two fundamental tools. The first one is a variational formulation of $\lambda_1(\Omega)$, which relies heavily on the symmetry of the Laplacian:
\begin{equation} \label{variational}
\lambda_1(\Omega)=\min_{u\in H^1_0(\Omega)\backslash\{0\}} \frac{\displaystyle{\int_{\Omega}} \left\vert \nabla u(x)\right\vert^2dx}{\displaystyle{\int_{\Omega}} u(x)^2dx}.
\end{equation}
The second ingredient of the proof is a spherical rearrangement argument, namely the Schwarz symmetrization of functions. Namely, if $u:\Omega\rightarrow \R$, the Schwarz spherical rearrangement of $u$ is the nonnegative function $u^{\ast}:\Omega^{\ast}\rightarrow \R$ which is radially non-increasing (which means, more precisely, that, if $R>0$ is the radius of $\Omega^{\ast}$, there exists a non-increasing function $v:\left[0,R\right]\rightarrow \R$ such that $u^{\ast}(x)=v(\left\vert x\right\vert)$ for all $x\in\Omega^{\ast}$) and satisfies
\[
\left\vert \left\{x\in \Omega;\ \left\vert u(x)\right\vert>\lambda\right\}\right\vert=\left\vert \left\{x\in \Omega^{\ast};\ u^{\ast}(x)>\lambda\right\}\right\vert
\]
for all $\lambda>0$. When $u\in H^1_0(\Omega)$, one has $u^{\ast}\in H^1_0(\Omega)$, $\left\Vert u\right\Vert_{L^2(\Omega)}=\left\Vert u^{\ast}\right\Vert_{L^2(\Omega)}$ and $\left\Vert \nabla u^{\ast}\right\Vert_{L^2(\Omega)}\leq \left\Vert \nabla u\right\Vert_{L^2(\Omega)}$ (see \cite{polyaszego}). Inequality (\ref{RFK}) follows immediately from (\ref{variational}) and these properties of $u^{\ast}$. 

When $v\neq 0$, the operator $-\Delta+v\cdot \nabla$ is non-self-adjoint, and there is no simple variational formulation of its first eigenvalue such as (\ref{variational}) --min-max formulations of the pointwise type (see \cite{bnv} and Section \ref{secminmax}) or of the integral type (see \cite{holland}) certainly hold, but they do not seem to help in our context. Therefore, it seems impossible to adapt the ``classical'' proof to prove Theorem \ref{th1}. However, the proof of Theorem \ref{th1} also relies on a new type of rearrangement argument, which is definitely different from the usual rearrangement of functions. 

First, using essentially the maximum principle and the Hopf lemma, one establishes Theorem \ref{th2} and Theorem \ref{th3}. Once this is done, we are reduced to proving that $\underline{\lambda}(\Omega,\tau)\geq \underline{\lambda}(\Omega^{\ast},\tau)$ for all $\tau\ge 0$, and that equality holds only when $\Omega$ is a ball.

To that purpose, we consider the function $\underline{\varphi}$ satisfying (\ref{eqmin}), and define a suitable rearrangement of $\underline{\varphi}$, which is different from the spherical symmetrization mentioned before. Let us briefly explain what the idea of this rearrangement is. Denote by $R$ the radius of $\Omega^{\ast}$. For all $0\leq a<1$, define
\[
\Omega_a=\left\{x\in \Omega,\ a<\underline{\varphi}(x)\le 1\right\}
\]
and define $\rho(a)\in\left(0,R\right]$ such that $\left\vert \Omega_a\right\vert=\left\vert B^n_{\rho(a)}\right\vert$. Define also $\rho(1)=0$. The function $\rho:\left[0,1\right]\rightarrow \left[0,R\right]$ is decreasing, continuous, one-to-one and onto. Then, the rearrangement of $\underline{\varphi}$ is the radially decreasing function $u:\Omega^{\ast}\rightarrow \R$ vanishing on $\partial\Omega^{\ast}$ such that, for all $0\leq a<1$,
\[
\int_{\Omega_a} \Delta \underline{\varphi}(x)dx=\int_{B^n_{\rho(a)}} 
\Delta u(x)dx.
\]
The fundamental inequality satisfied by $u$, which is also the key point in the proof of Theorem \ref{th1}, is the fact that, for all $x\in \overline{\Omega^{\ast}}$, 
\begin{equation} \label{fund}
u(x)\geq \rho^{-1}(\left\vert x\right\vert)
\end{equation}
(see Corollary \ref{cor36} below). Strictly speaking, this fact is not completely correct, because the function $\underline{\varphi}$ is not regular enough, and we have to deal with suitable approximations of $\underline{\varphi}$. We refer to Section \ref{sec3} for rigorous statements and proofs. Let us just mention that the proof of (\ref{fund}) relies, among other things, on the usual isoperimetric inequality in $\R^n$. From (\ref{fund}) and arguments involving the maximum principle and the Hopf lemma again, we derive the conclusion if Theorem \ref{th1}.

Finally, using again the same construction and the isoperimetric inequality, we prove that equality in Theorem \ref{th1} is attained if, and only if, up to translation, $\Omega=\Omega^{\ast}$ and $v=\tau e_r$.

\begin{rem}{\rm Notice that the proof of Theorem \ref{th1} given in Section \ref{sec3} still works for $\tau=0$ and then provides an alternative proof of the Rayleigh-Faber-Krahn inequality (\ref{RFK}) for the Dirichlet Laplacian.}
\end{rem}

\noindent{\bf{Outline of the paper.}} The paper is organized as follows: in Section \ref{sec2}, we show Theorem \ref{th2}, Theorem \ref{th3} and various properties of $\underline{\lambda}(\Omega,\tau)$ and $\overline{\lambda}(\Omega,\tau)$. Using the previous results and the rearrangement argument briefly described above, we prove Theorem \ref{th1} in Section \ref{sec3}. Finally, we prove in Appendix the results mentionned in Remark \ref{rem2} and Proposition \ref{prosup}.


\SE{Optimization problems in fixed domains}\label{sec2}


\subsection{Proof of Theorem \ref{th2}}

Throughout this section, we fix $\Omega\in {\mathcal C}$ of class $C^{2,\alpha}$ 
with $0<\alpha<1$ and $\tau\geq 
0$. The proof of Theorem \ref{th2} relies on the following comparison 
lemma:

\begin{lem} \label{comparison}
Let $\mu\in\R$ and $v\in L^{\infty}(\Omega,\R^n)$. Assume that 
$\varphi$ and $\psi$ are functions in 
$W^{2,p}(\Omega)$ for all $1\leq p<+\infty$, vanishing on 
$\partial\Omega$, satisfying $\left\Vert \varphi\right\Vert_{\infty}=
\left\Vert \psi\right\Vert_{\infty}$. Assume also that $\varphi\geq 
0$ in $\Omega$, $\psi>0$ in $\Omega$ and
\[
\left\{
\begin{array}{l}
-\Delta\psi+v.\nabla\psi\geq \mu\psi\mbox{ a. e. in }\Omega,\\
\\
-\Delta\varphi+v.\nabla \varphi\leq \mu\varphi\mbox{ a. e. in }\Omega.
\end{array}
\right.
\]
Then $\varphi=\psi$ in $\Omega$.
\end{lem}

\noindent{\bf Proof. }Since $\psi$ is not constant in $\Omega$, the Hopf 
lemma yields $\displaystyle \frac{\partial \psi}{\partial \nu}<0$ on 
$\partial\Omega$, where, for all $x\in \partial\Omega$, 
$\displaystyle \frac{\partial \psi}{\partial \nu}(x)=\nabla 
\psi(x)\cdot\nu(x)$ and $\nu(x)$ denotes the outward normal unit 
vector 
at $x$. Since $\varphi\in C^{1,\beta}(\overline{\Omega})$ for all 
$0\leq \beta<1$, $\varphi\geq 0$ in $\Omega$ and $\varphi=0$ on  
$\partial 
\Omega$, it follows that there exists $\gamma>0$ such that 
$\gamma\psi>\varphi$ in $\Omega$. Define
\[
\gamma^{\ast}=\inf\left\{\gamma>0,\ \gamma\psi>\varphi\mbox{ in 
}\Omega\right\}.
\]
One clearly has $\gamma^{\ast}\psi\geq \varphi$ in $\Omega$, so that 
$\gamma^{\ast}>0$. Define $w=\gamma^{\ast}\psi-\varphi$ and assume 
that 
$w>0$ everywhere in $\Omega$. Since 
\begin{equation} \label{diff}
-\Delta w+v\cdot \nabla w-\mu w\geq 0\mbox{ in }\Omega
\end{equation}
and $w=0$ on $\partial\Omega$, the Hopf maximum principle implies 
that $\displaystyle \frac{\partial w}{\partial \nu}<0$. As above, this yields the existence of $\kappa>0$ such that $w>\kappa\varphi$ in $\Omega$, whence
$$\frac{\gamma^{\ast}}{1+\kappa}\psi>\varphi\ \hbox{ in }\Omega.$$
This is a contradiction with the minimality of 
$\gamma^{\ast}$.

Thus, there exists $x_0\in \Omega$ such that $w(x_0)=0$ (i.e.  
$\gamma^{\ast}\psi(x_0)=\varphi(x_0)$). It follows from (\ref{diff}), 
the fact 
that $w\geq 0$ in $\Omega$, $w(x_0)=0$ and from the strong maximum 
principle, that $w=0$ in 
$\Omega$, which means that $\varphi$ and $\psi$ are proportional. 
Since 
they are non-negative in $\Omega$ and have the same $L^{\infty}$ norm 
in $\Omega$, 
one has $\varphi=\psi$, which ends the proof of Lemma 
\ref{comparison}. \hfill\fin\break

We now turn to the proof of assertion $(a)$ in Theorem \ref{th2} and 
begin 
with the existence of $\underline{v}$:

\begin{lem} \label{existence}
There exists $\underline{v}\in L^{\infty}(\Omega,\R^n)$ with 
$\left\Vert 
\underline{v}\right\Vert_{\infty}=\tau$ such that 
$\underline{\lambda}(\Omega,\tau)=\lambda_1(\Omega,\underline{v})$ $(>0)$. 
\end{lem}

\noindent{\bf Proof. } In this proof, we write $\lambda$ instead of 
$\underline{\lambda}(\Omega,\tau)$. Let $(v_k)_{k\geq 1}$ be a 
sequence of $L^{\infty}(\Omega,\R^n)$ functions with 
$\left\Vert v_k\right\Vert_{\infty}\leq \tau$ such that 
$\lambda_1(\Omega,v_k)$ converges to 
$\underline{\lambda}$, and, for all $k\geq 1$, set 
$\lambda_k=\lambda_1(\Omega,v_k)$ and 
$\varphi_k=\varphi_{\Omega,v_k}$. For 
all $k\geq 1$, one has
\[
-\Delta\varphi_{k}+v_{k}\cdot\nabla\varphi_{k}=\lambda_{k}\varphi_{k}\mbox{ a. e. in }\Omega.
\]
Since the $v_k$'s are uniformly bounded in $L^{\infty}(\Omega,\R^n)$, since 
the $\varphi_{k}$'s are uniformly bounded in $L^p(\Omega)$ 
(say, $1<p<+\infty$) and since the 
$\lambda_k$'s are bounded, the Agmon-Douglis-Nirenberg estimates (see 
\cite{adn,gt}) show that the $\varphi_k$'s are uniformly 
bounded in $W^{2,p}(\Omega)$ for all $1<p<+\infty$. 
Therefore, up to a subsequence, there exist $\varphi$ belonging to 
$W^{2,p}(\Omega)$ for all $1\le p<+\infty$ and to 
$C^{1,\beta}(\overline{\Omega})$ 
for all $0\leq \beta<1$, and a function $f\in L^{\infty}(\Omega)$ 
such that 
$\varphi_k$ converges to $\varphi$ weakly in 
$W^{2,p}(\Omega)$ and strongly in $C^{1,\beta}(\overline{\Omega})$, 
and 
$v_k\cdot\nabla\varphi_k$ converges to $f$ for the $\ast$-weak 
topology of $L^{\infty}(\Omega)$. Observe that $\varphi=0$ on 
$\partial\Omega$, $\varphi\geq 0$ in $\Omega$ and $\left\Vert 
\varphi\right\Vert_{\infty}=1$. Since, for all $k\geq 1$, 
$-\Delta\varphi_k\leq \lambda_k\varphi_k+\tau\left\vert \nabla 
\varphi_k\right\vert$ a. e. in $\Omega$, one has
\[
\left\{
\begin{array}{l}
-\Delta\varphi\leq \lambda\varphi+\tau\left\vert 
\nabla\varphi\right\vert\mbox{ a. e. in }\Omega,\\
\\
-\Delta\varphi + f=\lambda\varphi\mbox{ a. e. in }\Omega.
\end{array}
\right.
\]
It follows that $f\geq -\tau\left\vert \nabla \varphi\right\vert$ a. 
e. in 
$\Omega$. For all $x\in\Omega$, define
\[
\underline{v}(x)=
\left\{
\begin{array}{ll}
\displaystyle -\tau\frac{\nabla\varphi(x)}{\left\vert \nabla 
\varphi(x)\right\vert} & \mbox{if }\nabla\varphi(x)\neq 0,\\
\\
0 & \mbox{if }\nabla\varphi(x)=0.
\end{array}
\right.
\]
Since $\varphi$ is not constant, one has $\left\Vert 
\underline{v}\right\Vert_{\infty}=\tau$. Set 
$\mu=\lambda_1(\Omega,\underline{v})$ and 
$\psi=\varphi_{\Omega,\underline{v}}$ (recall that $\psi\in 
W^{2,p}(\Omega)$ for all $1\leq p<+\infty$). The very definition of 
$\underline{\lambda}(\Omega,\tau)$ yields that 
$\lambda\leq \mu$. Moreover,
\[
-\Delta\varphi+\underline{v}\cdot\nabla \varphi=-\Delta\varphi-\tau\left\vert 
\nabla \varphi\right\vert\leq \lambda\varphi\leq \mu\varphi\mbox{ a. 
e. in }\Omega,
\]
whereas $-\Delta\psi+\underline{v}\cdot\nabla \psi=\mu\psi$ a. e. in $\Omega$, 
$\psi$ is positive in $\Omega$, $\varphi$ is non-negative in 
$\Omega$, $\varphi=\psi=0$ on $\partial\Omega$ and $\left\Vert 
\varphi\right\Vert_{\infty}=\left\Vert \psi\right\Vert_{\infty}=1$. 
Lemma \ref{comparison} therefore yields $\varphi=\psi$ in $\Omega$. 
Thus,
\[
\mu\varphi=-\Delta\varphi+\underline{v}\cdot\nabla \varphi\leq \lambda\varphi\leq 
\mu\varphi\mbox{ a. e. in }\Omega,
\]
which means that $\lambda=\mu$, or, in other words, that 
$\underline{\lambda}(\Omega,\tau)=\lambda_1(\Omega,\underline{v})$. Lastly, $\lambda_1(\Omega,\underline{v})>0$ as already underlined in Section \ref{intro}.\hfill\fin\break

To prove the ``uniqueness'' statement in Theorem \ref{th2}, we first 
establish the following result:

\begin{lem} \label{uniqueness}
Let $v_0\in L^{\infty}(\Omega,\R^n)$ with $\left\Vert 
v_0\right\Vert_{\infty}\leq \tau$ such that 
$\lambda_1(\Omega,v_0)=\underline{\lambda}(\Omega,\tau)$ and 
call $\varphi=\varphi_{\Omega,v_0}$. Then $\varphi\in C^{2,\alpha}(\overline{\Omega})$ (up to the choice of the continuous representant in the class of $\varphi$), $\left\Vert 
v_{0}\right\Vert_{\infty}=\tau$ and 
$v_0\cdot\nabla\varphi=-\tau\left\vert \nabla \varphi\right\vert$ a. e. in $\Omega$.
\end{lem}

\noindent{\bf Proof.} Define, for all $x\in \Omega$,
\[
v(x)=
\left\{
\begin{array}{ll}
\displaystyle -\tau\frac{\nabla\varphi(x)}{\left\vert \nabla 
\varphi(x)\right\vert} & \mbox{if }\nabla\varphi(x)\neq 0,\\
\\
0 & \mbox{if }\nabla\varphi(x)=0,
\end{array}
\right.
\]
so that $\left\Vert v\right\Vert_{\infty}=\tau$ (indeed, $\varphi$ is 
not constant in $\Omega$). If 
$\lambda=\lambda_1(\Omega,v)$, one has $\lambda\geq 
\underline{\lambda}(\Omega,\tau)$. Moreover, for 
$\psi=\varphi_{\Omega,v}$, one has
\[
\left\{
\begin{array}{l}
-\Delta\varphi+v\cdot\nabla\varphi=-\Delta\varphi-\tau\left\vert 
\nabla 
\varphi\right\vert\leq 
-\Delta\varphi+v_0\cdot\nabla\varphi=\underline{\lambda}(\Omega,\tau)\varphi\leq 
\lambda \varphi \mbox{ a. e. in }\Omega,\\
\\
-\Delta\psi+v\cdot\nabla \psi=\lambda\psi\mbox{ a. e. in }\Omega.
\end{array}
\right.
\]
Since $\psi>0$ in $\Omega$, $\varphi>0$ in $\Omega$, $\psi=\varphi=0$ 
on $\partial\Omega$ and $\left\Vert 
\psi\right\Vert_{\infty}=\left\Vert \varphi\right\Vert_{\infty}=1$, 
Lemma 
\ref{comparison} ensures that $\psi=\varphi$. As a consequence, 
$\lambda=\underline{\lambda}(\Omega,\tau)$, and $v_0\cdot\nabla 
\varphi=-\tau\left\vert \nabla \varphi\right\vert$ a. e. in $\Omega$. Therefore, since $\|v_0\|_{\infty}\le\tau$ and $\varphi$ is not constant, one gets $\|v_0\|_{\infty}=\tau$. Furthermore, up to the choice of the continuous representant in the class of $\varphi$, $\varphi$ satisfies
$$-\Delta\varphi=\tau|\nabla\varphi|+\underline{\lambda}(\Omega,\tau)\varphi\in C^{0,\alpha}(\overline{\Omega}),$$
hence $\varphi\in C^{2,\alpha}(\overline{\Omega})$ from Schauder estimates \cite{gt}. That ends the proof of Lemma~\ref{uniqueness}.\hfill\fin\break
    
The last result for the proof of Theorem \ref{th2} is the following 
one:

\begin{lem} \label{uniquenessbis}
Let $\lambda\in\R$ and $\varphi\in C^{2,\alpha}(\overline{\Omega})$, such that $\varphi=0$ on $\partial\Omega$, $\varphi>0$ in 
$\Omega$, $\left\Vert \varphi\right\Vert_{\infty}=1$ and
\[
-\Delta\varphi-\tau\left\vert 
\nabla\varphi\right\vert=\lambda\varphi\mbox{ in }\Omega.
\]
Then $\lambda=\underline{\lambda}(\Omega,\tau)$ and, if $v\in 
L^{\infty}(\Omega,\R^n)$ is such that $\left\Vert 
v\right\Vert_{\infty}=\tau$ and 
$\underline{\lambda}(\Omega,\tau)=\lambda_1(\Omega,v)$, then 
$\varphi=\varphi_{\Omega,v}$.
\end{lem}

\noindent{\bf Proof.} Let $v\in L^{\infty}(\Omega,\R^n)$ such that 
$\left\Vert v\right\Vert_{\infty}=\tau$ and 
$\underline{\lambda}(\Omega,\tau)=\lambda_1(\Omega,v)$ (such a $v$ exists thanks to Lemma \ref{existence}), and set 
$\psi=\varphi_{\Omega,v}$, so that
\[
-\Delta\psi+v\cdot \nabla\psi=-\Delta\psi-\tau\left\vert 
\nabla\psi\right\vert=\underline{\lambda}(\Omega,\tau)\psi\mbox{ in }\Omega
\]
from Lemma \ref{uniqueness}. Define also
\[
w(x)=
\left\{
\begin{array}{ll}
\displaystyle -\tau\frac{\nabla\varphi(x)}{\left\vert \nabla 
\varphi(x)\right\vert} & \mbox{if }\nabla\varphi(x)\neq 0,\\
\\
0 & \mbox{if }\nabla\varphi(x)=0.
\end{array}
\right.
\]
One has $\left\Vert w\right\Vert_{\infty}=\tau$ and
\[
-\Delta\varphi+w\cdot \nabla\varphi=-\Delta \varphi-\tau\left\vert 
\nabla \varphi\right\vert=\lambda\varphi\mbox{ in }\Omega,
\]
so that, by uniqueness of the first eigenvalue and eigenfunction for $-\Delta +w\cdot\nabla$, one has $\varphi=\varphi_{\Omega,w}$ and $\lambda=\lambda_1(\Omega,w)\geq\underline{\lambda}(\Omega,\tau)$. As a consequence,
\[
-\Delta\varphi+v\cdot \nabla\varphi\geq -\Delta\varphi-\tau\left\vert 
\nabla\varphi\right\vert=\lambda\varphi\geq 
\underline{\lambda}(\Omega,\tau)\varphi\mbox{ in }\Omega.
\]
Another application of Lemma \ref{comparison} shows that 
$\varphi=\psi=\varphi_{\Omega,v}$, and that 
$\lambda=\underline{\lambda}(\Omega,\tau)$.\hfill\fin\break

We now complete the\hfill\break
{\bf{Proof of Theorem \ref{th2}.}} The existence of 
$\underline{v}$ and the identity $\underline{v}\cdot 
\underline{\varphi}=-\tau\left\vert\nabla\underline{\varphi}\right\vert$ in $\Omega$ have already been proved. 
Let $v_1$ and $v_2\in L^{\infty}(\Omega)$ with $\left\Vert 
v_1\right\Vert_{\infty}\leq \tau$, $\left\Vert 
v_2\right\Vert_{\infty}\leq \tau$ and 
$\lambda_1(\Omega,v_1)=\lambda_1(\Omega,v_2)=\underline{\lambda}(\Omega,\tau)$. 
Note $\varphi_1=\varphi_{\Omega,v_1}$ and 
$\varphi_2=\varphi_{\Omega,v_2}$, so that $\varphi_1$, $\varphi_2\in C^{2,\alpha}(\overline{\Omega})$, 
$v_1\cdot\nabla\varphi_1=-\tau\left\vert \nabla\varphi_1\right\vert$ 
and $v_2\cdot\nabla\varphi_2=-\tau\left\vert 
\nabla\varphi_2\right\vert$ a. e. in $\Omega$ by Lemma 
\ref{uniqueness}. Furthermore, $\|v_1\|_{\infty}=\|v_2\|_{\infty}=\tau$. One has
\be\label{v12}\left\{\begin{array}{l}
-\Delta\varphi_{1}+v_1\cdot\nabla\varphi_{1}=-\Delta\varphi_1-\tau\left\vert \nabla
\varphi_1\right\vert=\underline{\lambda}(\Omega,\tau)\varphi_{1}\mbox{ in }\Omega,\\
\\
-\Delta\varphi_{2}+v_2\cdot\nabla\varphi_{2}= 
-\Delta\varphi_{2}-\tau\left\vert 
\nabla\varphi_{2}\right\vert=\underline{\lambda}(\Omega,\tau)\varphi_{2}\mbox{ in }\Omega.
\end{array}
\right.
\ee
Lemma \ref{uniquenessbis} therefore shows that 
$\varphi_{1}=\varphi_{2}:=\varphi$, which yields $v_1(x)=v_2(x)=-\tau\nabla\varphi(x)/|\nabla\varphi(x)|$ almost everywhere in the set of $x\in \Omega$ such that $\nabla \varphi(x)\neq 0$. What remains 
to be proved is the fact that $\nabla\varphi(x)\neq 0$ for almost every 
$x\in \Omega$ (that will then imply the uniqueness of $\underline{v}$ in Theorem \ref{th2} and the fact that $|\underline{v}(x)|=\tau$ almost everywhere in $\Omega$). To that purpose, we recall that, for all $p\ge 1$ 
and all function $g\in W^{1,p}(\Omega)$, $\nabla g=0$ almost 
everywhere in the set $\left\{x\in\Omega,\ g(x)=0\right\}$. For each 
$1\leq i\leq n$, this observation applied with $\displaystyle 
\frac{\partial \varphi}{\partial x_i}$ yields
\[
\int_{\left\{\nabla \varphi=0\right\}}\Delta 
\varphi(x)dx=\sum\limits_{i=1}^n \int_{\Omega} 
\frac{\partial^2\varphi}{\partial x_i^2}(x){\bf 
1}_{\left\{\nabla\varphi=0\right\}}(x)dx=0,
\]
where ${\bf{1}}_E$ denotes the characteristic function of a set $E$. Therefore, since $\underline{\lambda}(\Omega,\tau)>0$, (\ref{v12}) 
ensures that
\[
\int_{\left\{\nabla \varphi=0\right\}} \varphi(x)dx=0,
\]
and since $\varphi(x)>0$ for all $x\in \Omega$, one has $\left\vert 
\left\{\nabla\varphi=0\right\}\right\vert=0$, which ends the proof of 
assertion $(a)$ in Theorem \ref{th2}. 

The proof of assertion $(b)$ is 
entirely similar, except for the proof of the existence of $\overline{v}$, for which we need to know a priori that, given
$\Omega\in {\mathcal C}$ and $\tau\geq 0$, there exists $C>0$ such that, for all $v\in L^{\infty}(\Omega,\R^n)$ satisfying
$\left\Vert v\right\Vert_{\infty}\leq \tau$, one has $\lambda_1(\Omega,v)\leq C$. But this is true in virtue of Proposition 5.1
in \cite{bnv}, which yields the existence of some constant $C(\Omega,\tau)\ge 0$ such that $\left\vert \lambda_1(\Omega,v)-\lambda_1(\Omega,0)\right\vert \leq C(\Omega,\tau)\left\Vert
v\right\Vert_{\infty}$.\hfill\fin


\subsection{The case of a ball}

This section is devoted to the\hfill\break
\noindent{\bf{Proof of Theorem \ref{th3}.}} Denote by $R>0$ the radius of $\Omega$, namely $\Omega=B^n_R$, and let $\varphi=\varphi_{\Omega,\underline{v}}$. Recall that
\[
-\Delta \varphi-\tau\left\vert 
\nabla \varphi\right\vert=\underline{\lambda}(\Omega,\tau)\varphi\mbox{ in }\Omega,
\]
thus, $\varphi\in C^{2,\beta}(\overline{\Omega})$ for all $0\le\beta<1$ from Schauder estimates \cite{gt}. Let $A$ be an orthogonal transformation in $\R^n$ and 
$\psi=\varphi\circ A$. Easy computations show that $\psi$ satisfies 
the same equation as $\varphi$, and Lemma \ref{uniquenessbis} ensures that 
$\psi=\varphi$ in $\Omega$. In other words, $\varphi$ is radial, and 
we may define $u:\left[0,R\right]\rightarrow \R$ such that, for all 
$x\in \Omega$, $\varphi(x)=u(\left\vert x\right\vert)$. The function $u$ is continuous in $\left[0,R\right]$, positive in
$\left[0,R\right)$, its maximum in $\left[0,R\right]$ is $1$ and $u(R)=0$.

We claim that $u$ is decreasing in $\left[0,R\right]$.  First, there exists $r_0\in \left[0,R\right)$ such that $u(r_0)=1$. If
$r_0>0$, since 
\begin{equation} \label{eigenvalue}
-\Delta\varphi+\underline{v}\cdot \nabla\varphi=\underline{\lambda}(\Omega,\tau)\varphi>0
\end{equation}
in $B^n_{r_0}$, the maximum principle shows that $\varphi(x)\geq 1$ for all $x\in B^n_{r_0}$, and therefore that $\varphi(x)=1$ in
$B^n_{r_0}$, which contradicts (\ref{eigenvalue}). Thus, $u(0)=\varphi(0)=1$ and $u(r)<1$ for all $0<r\leq R$. Assume now that $u$
is not decreasing in $\left[0,R\right]$. This means that there exist $0\leq r_0<r_1\leq R$ such that $u(r_0)\leq u(r_1)$. Notice
that $r_1<R$, $0<r_0$ and set $m=u(r_1)\in \left(0,1\right)$. Since $u$ is continuous, there exists $r_2\in \left(0,r_0\right]$
such that $u(r_2)=m$ and $m<u(r)\le 1$ for all $r\in[0,r_2)$. Since (\ref{eigenvalue}) holds in the spherical shell $U=\left\{x\in \R^n;\ r_2<\left\vert
x\right\vert<r_1\right\}$, the maximum principle applied to $\varphi$ in $U$  shows that $\varphi(x)\geq m$ for all $x\in U$. If $\varphi(x_0)=m$ for some $x_0\in U$, the the strong maximum principle implies that $\varphi$ is constant in $U$ and this is impossible because of (\ref{eigenvalue}). Therefore, $\varphi(x)>m$ in $U$, and Hopf lemma yields
$$u'(r_2)=\nabla\varphi(x)\cdot e_r(x)>0$$
for all $x$ with $|x|=r_2$. This contradicts the definition of $r_2$.

Finally, $u$ is decreasing in $\left[0,R\right]$ and Hopf lemma applied to (\ref{eigenvalue}) on the boundary of any ball $B^n_r$ with $r\in(0,R]$ implies that $\nabla\varphi(x)\cdot e_r(x)<0$ for all $x\in\overline{\Omega}\backslash\{0\}$. To sum up, $\nabla\varphi(x)=u^{\prime}(\left\vert x\right\vert)e_r\neq 0$ for all $x\in 
\overline{\Omega}\setminus \left\{0\right\}$ and, for all $x\neq 0$, 
\[
\underline{v}(x)=-\tau \frac{\nabla\varphi(x)}{\left\vert 
\nabla\varphi(x)\right\vert}=\tau e_r(x)
\]
from Lemma \ref{uniqueness}.

The argument for $\overline{v}$ is completely similar, and this ends the proof of Theorem \ref{th3}. \hfill\fin


\subsection{Further properties of $\underline{\lambda}(\Omega,\tau)$, $\overline{\lambda}(\Omega,\tau)$ and $F_n(m,\tau)$}\label{secfurther}
    
We prove here the continuity and monotonicity of the maps $\tau\mapsto\underline{\lambda}(\Omega,\tau)$ and $\tau\mapsto\overline{\lambda}(\Omega,\tau)$.

\begin{lem} \label{continuity}
Let $\Omega\in {\mathcal C}$ be fixed. Then the function 
$\tau\mapsto \underline{\lambda}(\Omega,\tau)$ is continuous and 
decreasing in $\left[0,+\infty\right)$.
\end{lem}

\noindent{\bf Proof.} That this function is non-increasing follows at once 
from the definition of $\underline{\lambda}(\Omega,\tau)$. Let now $0\le\tau<\tau'$ be given. Under the notations of part $(a)$ of Theorem \ref{th2}, there is a (unique) $\underline{v}\in L^{\infty}(\Omega,\R^n)$ such that $\underline{\lambda}(\Omega,\tau)=\lambda_1(\Omega,\underline{v})$ and $\|\underline{v}\|_{\infty}\le\tau$ (actually, $\|\underline{v}\|_{\infty}=\tau$). Since $\|\underline{v}\|_{\infty}<\tau'$, the uniqueness result in Theorem \ref{th2} yields $\lambda_1(\Omega,\underline{v})>\underline{\lambda}(\Omega,\tau')$. Therefore, the map $\tau\mapsto\underline{\lambda}(\Omega,\tau)$ in decreasing in $[0,+\infty)$.

Fix now $\tau\geq 0$, let $(\tau_k)_{k\geq 1}$ be a sequence of non-negative 
numbers converging to $\tau$ and write 
$\lambda_k=\underline{\lambda}(\Omega,\tau_k)$. For each $k\geq 1$, 
Theorem \ref{th2} provides the existence of a function 
$\varphi_k\in C^{2,\alpha}(\overline{\Omega})$ ($0<\alpha<1$ is so that $\Omega$ is of class $C^{2,\alpha}$) such that
\be\label{lambdak}
-\Delta\varphi_k-\tau_k\left\vert \nabla 
\varphi_k\right\vert=\lambda_k\varphi_k\mbox{ in }\Omega
\ee
with $\varphi_k>0$ in $\Omega$, $\varphi_k=0$ on $\partial\Omega$ and $\|\varphi_k\|_{\infty}=1$. 
Observe that the sequence $(\lambda_k)_{k\geq 1}$ is bounded (if $A>0$ is such that $0\leq \tau_k\leq A$ for all $k$, one has
$0<\underline{\lambda}(\Omega,A)\leq \lambda_k\leq \underline{\lambda}(\Omega,0)$) and therefore 
the $\varphi_k$'s are uniformly bounded in $W^{2,p}(\Omega)$ for all $1\le p<+\infty$ and then, again because of (\ref{lambdak}), the functions $\varphi_k$ are uniformly bounded in $C^{2,\alpha}(\overline{\Omega})$. Up to a 
subsequence, the sequence $(\lambda_k)_{k\geq 1}$ converges to 
$\lambda>0$ and there exists a function $\varphi\in C^{2,\alpha}(\overline{\Omega})$, such that 
$(\varphi_k)_{k\geq 1}$ converges to $\varphi$ strongly in 
$C^{2,\beta}(\overline{\Omega})$ for all $0\leq \beta<\alpha$. One therefore has
\[
-\Delta\varphi-\tau\left\vert \nabla 
\varphi\right\vert=\lambda\varphi\mbox{ in }\Omega,
\]
and $\varphi\ge 0$ in $\overline{\Omega}$, $\|\varphi\|_{\infty}=1$, $\varphi=0$ on $\partial\Omega$. Since $-\Delta\varphi\ge 0$ in $\Omega$, the strong maximum principle yields $\varphi>0$ in $\Omega$. It follows from Lemma \ref{uniquenessbis} that 
$\lambda=\underline{\lambda}(\Omega,\tau)$. Thus, the sequence 
$(\lambda_k)_{k\geq 1}$ is bounded and any converging subsequence of 
$(\lambda_k)_{k\geq 1}$ converges to 
$\underline{\lambda}(\Omega,\tau)$, which shows that 
$\lambda_k\rightarrow \underline{\lambda}(\Omega,\tau)$ and ends the 
proof of Lemma \ref{continuity}. \hfill\fin\break

Similarly, the following result holds

\begin{lem} \label{continuitybis}
Let $\Omega\in {\mathcal C}$ be fixed. Then the function 
$\tau\mapsto \overline{\lambda}(\Omega,\tau)$ is continuous and 
increasing in $\left[0,+\infty\right)$.
\end{lem}
The proof is similar to the one of Lemma \ref{continuity}, except that 
the function $\varphi$ obtained in the end of the argument satisfies $-\Delta\varphi+\tau\left\vert 
\nabla\varphi\right\vert=\lambda\varphi$ in $\Omega$. Setting
\[
v(x)=\left\{
\begin{array}{ll}
\displaystyle \tau \frac{\nabla\varphi(x)}{\left\vert 
\nabla\varphi(x)\right\vert} &\mbox{ if }\nabla\varphi(x)\neq 0,\\
\\
0 & \mbox{ if }\nabla\varphi(x)=0,
\end{array}
\right.
\]
one has $-\Delta\varphi+v\cdot \nabla\varphi\geq 0$ on $\Omega$, and the 
strong maximum principle yields $\varphi>0$ on $\Omega$, and one 
concludes as in the proof of Lemma \ref{continuity}. \hfill\fin\break

\begin{rem}{\rm An immediate application of Proposition 5.1 in 
\cite{bnv} yields that, for given $\Omega\in{\mathcal{C}}$, the map $v\mapsto\lambda_1(\Omega,v)$ is continuous as well 
(with respect to the $L^{\infty}$ norm for $v$).}
\end{rem}

\begin{rem}\label{contm}{\rm From Theorem \ref{th3} and Lemma \ref{continuity}, one gets that $F_n(m,\tau)=\lambda_1(B^n_{(m/\alpha_n)^{1/n}},\tau e_r)$ is 
decreasing in $\tau$. Furthermore, because of the strict monotonicity of the first eigenvalue with respect to the inclusion of the domains, $F_n(m,\tau)$ is 
decreasing in $m$ as well.\par
Actually, the map $(m,\tau)\mapsto F_{n}(m,\tau)$ is continuous on 
$\left(0,+\infty\right)\times \left[0,+\infty\right)$. The proof of 
this fact is very similar to the one of Lemma \ref{continuity}.}
\end{rem}


\SE{Proof of the main Faber-Krahn type inequality}\label{sec3}

This section is devoted to the proof of Theorem \ref{th1}. For the sake of clarity, we divide the proof into two parts~: first, in Section \ref{sec31}, we prove that, for all $\tau\ge 0$ and $\Omega\in{\mathcal{C}}$, there holds $\lambda_1(\Omega,v)\ge\lambda_1(\Omega^*,\tau e_r)$ for all $v\in L^{\infty}(\Omega)$ with $\|v\|_{L^{\infty}(\Omega)}\le\tau$. One recalls that $\Omega^*$ denotes the open Euclidean ball of $\R^n$ with center $0$ and such that $|\Omega^*|=|\Omega|$. Then, in Section \ref{sec32}, we discuss the case of equality.


\subsection{Proof of inequality $\lambda_1(\Omega,v)\ge\lambda_1(\Omega^*,\tau e_r)$}\label{sec31}

Let $\tau\ge 0$ and $\Omega\in{\mathcal{C}}$ be fixed, of class $C^{2,\alpha}$ for some $0<\alpha<1$. Let $R>0$ be the radius of the ball $\Omega^*$, namely $\Omega^*=B^n_R$ and $R=(|\Omega|/\alpha_n)^{1/n}$. Recall that $e_r$ denotes the unit radial vector field: $e_r(x)=x/|x|$ for $x\neq 0$.

In order to prove the first part of Theorem \ref{th1}, one shall show that $\lambda_1(\Omega,v)\ge\lambda_1(\Omega^*,\tau e_r)$ for all $v\in L^{\infty}(\Omega)$ such that $\|v\|_{L^{\infty}(\Omega)}\le\tau$. Owing to the definition of $\underline{\lambda}(\Omega,\tau)$, it is then enough to prove that
$$\underline{\lambda}(\Omega,\tau)\ge\lambda_1(\Omega^*,\tau e_r).$$

From the characterization of $\underline{\lambda}(\Omega,\tau)$ in Theorem \ref{th2}, let us call $\varphi$ the unique solution of (\ref{eqmin}), such that $\varphi>0$ in $\Omega$, $\|\varphi\|_{L^{\infty}(\Omega)}=1$ and $\varphi=0$ on $\partial\Omega$. Namely, the function $\varphi$, of class $C^{2,\alpha}(\overline{\Omega})$, satisfies
\begin{equation}\label{eqminbis}\left\{\baa{rcll}
-\Delta\varphi & = & \tau|\nabla\varphi|+\underline{\lambda}(\Omega,\tau)\varphi=:f & \hbox{in }\Omega\\
\varphi & > & 0 & \hbox{in }\Omega\\
\varphi & = & 0 & \hbox{on }\partial\Omega.\eaa\right.
\end{equation}
Furthermore, as already underlined in Section \ref{intro}, $\underline{\lambda}(\Omega,\tau)$ is positive and, on $\partial\Omega=\{x\in\overline{\Omega}$, $\varphi(x)=0\}$, $\nabla\varphi\neq 0$ from Hopf Lemma. Therefore, the continuous function $f$ is positive in $\overline{\Omega}$ if $\tau>0$. If $\tau=0$, then $f$ is continuous in $\overline{\Omega}$, positive in $\Omega$, and it vanishes on $\partial\Omega$.

From Stone-Weierstrass Theorem, let us choose a sequence $(f_k)_{k\in\N}$ of polynomial functions such that
$$f_k\ \to\ f\hbox{ in }C^0(\overline{\Omega})\hbox{ as }k\to+\infty.$$
If $\tau>0$, one can assume without loss of generality that $f_k>0$ in $\overline{\Omega}$ for all $k\in\N$. If $\tau=0$, setting $\tilde{f}_k=f_k-\min_{\overline{\Omega}}f_k+1/(k+1)$ and renaming $\tilde{f}_k$ as $f_k$, one still has that $f_k\to f$ uniformly in $\overline{\Omega}$ and $f_k>0$ in $\overline{\Omega}$ for all $k\in\N$.\par
For each $k\in\N$, call $\varphi_k$ the unique solution in $C^{2,\alpha}(\overline{\Omega})$ of
$$\left\{\baa{rcll}
-\Delta\varphi_k & = & f_k & \hbox{in }\Omega\\
\varphi_k & = & 0 & \hbox{on }\partial\Omega.\eaa\right.$$
The maximum principle implies that $\varphi_k>0$ in $\Omega$ for all $k\in\N$.\par
We first work with a fixed $k$ (large enough if necessary) and we will then pass to the limit as $k\to+\infty$ at the end of the proof.\par
Let first $k\in\N$ be fixed and let us introduce a few notations which will be used throughout this section. Call
$$M_k=\mathop{\max}_{x\in\overline{\Omega}}\ \varphi_k(x)$$
and, for $a\in[0,M_k]$,
$$\Sigma_{k,a}=\{x\in\overline{\Omega},\ \varphi_k(x)=a\}.$$\par
The function $\varphi_k$ can be written as $\varphi_k=\varphi'_k+\varphi''_k$ where $\varphi'_k$ is a polynomial function such that $-\Delta\varphi'_k=f_k$, and $\varphi''_k$ is harmonic in $\Omega$. Therefore, $\varphi_k$ is analytic in $\Omega$. On the other hand, Hopf lemma implies that $\frac{\partial\varphi_k}{\partial\nu}<0$ on $\partial\Omega$ ($\nu$ is the outward unit normal on $\partial\Omega$), whence the set $\{x\in\overline{\Omega},\ \nabla\varphi_k(x)=0\}$ is included in some compact set $K_k\subset\Omega$. It then follows that the set
$$Z_k=\{a\in[0,M_k],\ \exists x\in\Sigma_{k,a},\ \nabla\varphi_k(x)=0\}$$
of the critical values of $\varphi_k$ is finite (\cite{soucek}) and can then be written as
$$Z_k=\{a_{k,1},\cdots,a_{k,m_k}\}$$
for some $m_k\in\N^*$. Observe also that $M_k\in Z_k$ and that $0\not\in Z_k$. One can then assume that $0<a_{k,1}<\cdots<a_{k,m_k}=M_k$.\par
The set $Y_k=[0,M_k]\backslash Z_k$ of the non critical values of $\varphi_k$ is open relatively to $[0,M_k]$ and can be written as
$$Y_k=[0,M_k]\backslash Z_k=[0,a_{k,1})\cup(a_{k,1},a_{k,2})\cup\cdots\cup(a_{k,m_k-1},M_k).$$
For all $a\in Y_k$, the hypersurface $\Sigma_{k,a}$ is of class $C^2$ and $|\nabla\varphi_k|$ does not vanish in $\Sigma_{k,a}$. Therefore, in $Y_k$, the functions
\begin{equation}\label{ghi}\left\{\baa{l}
g_k\ :\ a\mapsto\displaystyle{\int_{\Sigma_{k,a}}}|\nabla\varphi_k(y)|^{-1}d\sigma_{k,a}(y)\\
h_k\ :\ a\mapsto\displaystyle{\int_{\Sigma_{k,a}}}f_k(y)|\nabla\varphi_k(y)|^{-1}d\sigma_{k,a}(y)\\
i_k\ :\ a\mapsto\displaystyle{\int_{\Sigma_{k,a}}}d\sigma_{k,a}(y)\eaa\right.\end{equation}
are continuous in $Y_k$, where $d\sigma_{k,a}$ denotes the surface measure on $\Sigma_{k,a}$ for $a\in Y_k$.\par
For all $a\in[0,M_k)$, let
$$\Omega_{k,a}=\{x\in\Omega,\ a<\varphi_k(x)\le M_k\}$$
and $\rho_k(a)\in(0,R]$ be defined so that
$$|\Omega_{k,a}|=|B^n_{\rho_k(a)}|=\alpha_n\rho_k(a)^n,$$
where one recalls that $\alpha_n$ is the volume of the unit ball $B^n_1$. One extends the function $\rho_k$ at $M_k$ by $\rho_k(M_k)=0$.

\begin{lem}\label{lem31} The function $\rho_k$ is a continuous decreasing map from $[0,M_k]$ onto $[0,R]$.
\end{lem}

\noindent{\bf{Proof.}} The function $\rho_k\ :\ [0,M_k]\to[0,R]$ is clearly decreasing since
$$\left|\{x\in\Omega,\ a<\varphi_k(x)\le b\}\right|>0$$
for all $0\le a<b\le M_k$.\par
Furthermore, for all $a\in(0,M_k]$ and all $1\leq i\leq n$, since $\displaystyle \frac{\partial^2\varphi_k}{\partial x_i^2}=0$ almost everywhere on the set where $\displaystyle \frac{\partial \varphi_k}{\partial x_i}=0$ as already mentioned in the proof of Theorem \ref{th2}, one has
$$\baa{rl}
\left.\displaystyle{\int_{\Sigma_{k,a}}}\Delta\varphi_k(x)dx\ =\ \displaystyle{\mathop{\sum}_{i=1}^n}\right( & \displaystyle{\int_{\Omega}}\displaystyle{\frac{\partial^2\varphi_k}{\partial x_i^2}}(x)\times{\bf{1}}_{\{\frac{\partial\varphi_k}{\partial x_i}>0\}}(x)\times{\bf{1}}_{\{\varphi_k=a\}}(x)dx\\
& \left.+\displaystyle{\int_{\Omega}}\displaystyle{\frac{\partial^2\varphi_k}{\partial x_i^2}}(x)\times{\bf{1}}_{\{\frac{\partial\varphi_k}{\partial x_i}<0\}}(x)\times{\bf{1}}_{\{\varphi_k=a\}}(x)dx\right).\eaa$$
But, using the same observation again, ${\bf{1}}_{\{\frac{\partial\varphi_k}{\partial x_i}<0\}}(x)\times{\bf{1}}_{\{\varphi_k=a\}}(x)={\bf{1}}_{\{\frac{\partial\varphi_k}{\partial x_i}>0\}}(x)\times{\bf{1}}_{\{\varphi_k=a\}}(x)=0$ for almost every $x\in\Omega$. As a consequence,
$$\int_{\Sigma_{k,a}}\Delta\varphi_k(x)dx=0.$$
But $-\Delta\varphi_k=f_k$ and the continuous function $f_k$ is positive in $\overline{\Omega}$. One then gets that $|\Sigma_{k,a}|=0$ for all $a\in(0,M_k]$. Notice that $|\Sigma_{k,0}|=|\partial\Omega|=0$ as well. Lastly, $\rho_k(0)=R$ and $\rho_k(M_k)=0$.\par
Therefore, the function $\rho_k$ is continuous in $[0,M_k]$. As a conclusion, $\rho_k$ it is then a one-to-one and onto map from $[0,M_k]$ to $[0,R]$.\hfill\fin

\begin{lem}\label{lem32} The function $\rho_k$ is of class $C^1$ in $Y_k$ and
$$\forall a\in Y_k,\quad\rho_k'(a)=-(n\alpha_n\rho_k(a)^{n-1})^{-1}g_k(a),$$
where the function $g_k$ is defined in (\ref{ghi}).
\end{lem}

\noindent{\bf{Proof.}} Fix $a\in Y_k$. Let $\eta>0$ be such that $[a,a+\eta]\subset Y_k$. For $t\in(0,\eta)$,
$$\baa{rcl}
\alpha_n(\rho_k(a+t)^n-\rho_k(a)^n)=|\Omega_{k,a+t}|-|\Omega_{k,a}| & = & -\displaystyle{\int_{\{a<\varphi_k(x)\le a+t\}}}dx\\
& = & -\displaystyle{\int_a^{a+t}}\left(\displaystyle{\int_{\Sigma_{k,b}}}|\nabla\varphi_k(y)|^{-1}d\sigma_{k,b}(y)\right)db\eaa$$
from the co-area formula. Hence,
$$\frac{\alpha_n(\rho_k(a+t)^n-\rho_k(a)^n)}{t}\to -g_k(a)\ \hbox{ as }t\to 0^+$$
for all $a\in Y_k$, due to the continuity of $g_k$ in $Y_k$. Similarly, one has that
$$\frac{\alpha_n(\rho_k(a+t)^n-\rho_k(a)^n)}{t}\to -g_k(a)\ \hbox{ as }t\to 0^-$$
for all $a\in Y_k\backslash\{0\}$.\par
The conclusion of the lemma follows since $Y_k\subset[0,M_k)$, whence $\rho_k(a)\neq 0$ for all $a\in Y_k$.\hfill\fin\break

The key-point in the proof of Theorem \ref{th1} is the construction of some auxiliary functions $u_k$ defined in $\Omega^*$. For each $k\in\N$, the function $u_k$ is obtained from $\varphi_k$ by a special new type of symmetrization.\par
Remember that $\Omega^*=B^n_R$ and define first, for all $r\in(0,R]$,
$$v_k(r)=\frac{1}{n\alpha_nr^{n-1}}\int_{\Omega_{k,\rho_k^{-1}(r)}}\Delta\varphi_k(x)dx,$$
where $\rho_k^{-1}\ :\ [0,R]\to[0,M_k]$ denotes the reciprocal of the function $\rho_k$. Set $v_k(0)=0$.

\begin{lem}\label{lem33} The function $v_k$ is continuous in $[0,R]$, and negative in $(0,R]$.
\end{lem}

\noindent{\bf{Proof.}} The continuity of $v_k$ in $(0,R]$ is a consequence of Lemma \ref{lem31} and the fact that $\Delta\varphi_k$ is continuous and thus bounded in $\overline{\Omega}$.\par
For $0<r\le R$, one has
$$|v_k(r)|\le(n\alpha_nr^{n-1})^{-1}\|\Delta\varphi_k\|_{\infty}\
\alpha_nr^n=n^{-1}\|\Delta\varphi_k\|_{\infty}\ r,$$
thus $v_k$ is continuous at $0$ as well.\par
The negativity of $v_k$ in $(0,R]$ follows from the negativity of $\Delta\varphi_k$ (in other words the positivity of $f_k$) in $\overline{\Omega}$ and the fact that $|\Omega_{k,a}|>0$ for all $a\in[0,M_k)$.\hfill\fin\break

For all $x\in\overline{\Omega^*}$, set
$$u_k(x)=-\int_{|x|}^Rv_k(r)dr.$$
The function $u_k$ is then radially symmetric in $\Omega^*$, decreasing in the variable $|x|$, positive in $\Omega^*$, vanishing on $\partial\Omega^*$ and, from Lemma \ref{lem33} and the fact that $v_k(0)=0$, $u_k$ is of class $C^1(\overline{\Omega^*})$.\par
Call
$$E_k=\{x\in\overline{\Omega^*},\ |x|\in\rho_k(Y_k)\}.$$
The set $E_k$ is a finite union of spherical shells and from Lemma \ref{lem31}, it is open relatively to $\overline{\Omega^*}$ and can be written as
$$E_k=\{x\in\R^n,\ |x|\in(0,\rho_k(a_{k,m_k-1}))\cup\cdots\cup(\rho_k(a_{k,2}),\rho_k(a_{k,1}))\cup(\rho_k(a_{k,1}),R]\}.$$
with $0=\rho_k(a_{k,m_k})=\rho_k(M_k)<\rho_k(a_{k,m_k-1})<\cdots<\rho_k(a_{k,1})<R$. Notice that $0\not\in E_k$.

\begin{lem}\label{lem34} The function $u_k$ is of class $C^2$ in $E_k$.
\end{lem}

\noindent{\bf{Proof.}} By definition of $u_k$, it is enough to prove that the function
$$w_k\ :\ r\mapsto\int_{\Omega_{k,\rho_k^{-1}(r)}}\Delta\varphi_k(x)dx=-\int_{\Omega_{k,\rho_k^{-1}(r)}}f_k(x)dx$$
is of class $C^1$ in $\rho_k(Y_k)$ ($\subset(0,R]$).\par
Let $r$ be fixed in $\rho_k(Y_k)=(0,\rho_k(a_{k,m-1}))\cup\cdots\cup(\rho_k(a_{k,2}),\rho_k(a_{k,1}))\cup(\rho_k(a_{k,1}),R]$ and let $\eta>0$ be such that $[r-\eta,r]\subset\rho_k(Y_k)$. For $t\in(0,\eta)$, one has
$$\baa{rcl}
w_k(r-t)-w_k(r) & = & \displaystyle{\int_{\{\rho_k^{-1}(r)<\varphi_k(x)\le\rho_k^{-1}(r-t)\}}}f_k(x)dx\\
& = & \displaystyle{\int_{\rho_k^{-1}(r)}^{\rho_k^{-1}(r-t)}}\left(\displaystyle{\int_{\Sigma_{k,a}}}f_k(y)|\nabla\varphi_k(y)|^{-1}d\sigma_{k,a}(y)\right)da=\displaystyle{\int_{\rho_k^{-1}(r)}^{\rho_k^{-1}(r-t)}}h_k(a)da,\eaa$$
where $h_k$ is defined in (\ref{ghi}). Since $\rho_k^{-1}$ is of class $C^1$ in $\rho_k(Y_k)$ from Lemma \ref{lem32} and since $h_k$ is continuous in $Y_k$, it follows that
$$\frac{w_k(r-t)-w_k(r)}{-t}\to h_k(\rho_k^{-1}(r))(\rho_k^{-1})'(r)=-\frac{n\alpha_nr^{n-1}h_k(\rho_k^{-1}(r))}{g_k(\rho_k^{-1}(r))}\hbox{ as }t\to 0^+.$$
The same limit holds as $t\to 0^-$ for all $r\in\rho_k(Y_k)\backslash\{R\}$. Therefore, the function $w_k$ is differentiable in $\rho_k(Y_k)$ and
$$w_k'(r)=-\frac{n\alpha_nr^{n-1}h_k(\rho_k^{-1}(r))}{g_k(\rho_k^{-1}(r))}\ \hbox{ for all }r\in\rho_k(Y_k).$$
Since $\rho_k^{-1}$ is continuous in $[0,R]$, and $g_k$ and $h_k$ are continuous in $Y_k$, the function $w_k$ is of class $C^1$ in $\rho_k(Y_k)$. That completes the proof of Lemma \ref{lem34}.\hfill\fin\break

The central argument is the following pointwise inequality satisfied by the ``symmetrized'' function $u_k$ obtained from $\varphi_k$. This inequality has its own independent interest. Besides the definition of $u_k$, it uses the classical isoperimetric inequality.

\begin{pro}\label{pro35} For all $x\in E_k$ and for all $\omega\ge 0$,
\begin{equation}\label{eq1pro35}
\Delta u_k(x)+\omega\ |\nabla u_k(x)|\ge\mathop{\min}_{y\in\Sigma_{k,\rho_k^{-1}(|x|)}}(\Delta\varphi_k(y)+\omega\
|\nabla\varphi_k(y)|).
\end{equation}
Furthermore, for any unit vector $e$ of $\R^n$, the function
$$\baa{rcl}
U_k\ :\ [0,M_k] & \to & \R_+\\
a & \mapsto & u_k(\rho_k(a)e)\eaa$$
is continuous in $[0,M_k]$, differentiable in $Y_k$, and
\begin{equation}\label{eq2pro35}
\forall a\in Y_k,\quad U_k'(a)\ge 1.
\end{equation}
\end{pro}

The first part of this proposition means that, for each $x\in E_k$ and $\omega\ge 0$, there is a point $y\in\overline{\Omega}$ such that $\varphi_k(y)=\rho_k^{-1}(|x|)$ and
$$\Delta u_k(x)+\omega\ |\nabla u_k(x)|\ge\Delta\varphi_k(y)+\omega\ |\nabla\varphi_k(y)|.$$\par
Before doing the proof of Proposition \ref{pro35}, let us first state the following

\begin{cor}\label{cor36} For all $x\in\overline{\Omega^*}$,
$$u_k(x)\ge\rho_k^{-1}(|x|).$$
\end{cor}

\noindent{\bf{Proof.}} On the one hand, $Y_k=[0,a_{k,1})\cup(a_{k,1},a_{k,2})\cup\cdots\cup(a_{k,m_k-1},M_k)$ and the function $\rho_k$ is continuous in $[0,M_k]$, differentiable in $Y_k$ (Lemmata \ref{lem31} and \ref{lem32}) and $\rho_k(0)=R$. On the other hand, the function $u_k$ is of class $C^1(\overline{\Omega^*})$, radially symmetric, $u_k(x)=0$ as soon as $|x|=R$. Corollary \ref{cor36} then follows from (\ref{eq2pro35}) and mean value theorem.\hfill\fin\break

\noindent{\bf{Proof of Proposition \ref{pro35}.}} Fix $x\in E_k$ and $\eta>0$ such that $\{y,\ |x|-\eta\le|y|\le|x|\}\subset E_k$. Call $r=|x|$. One has
$$\int_{\Sigma_{k,\rho_k^{-1}(r)}}|\nabla\varphi_k(y)|d\sigma_{k,\rho_k^{-1}(r)}(y)=-\int_{\partial\Omega_{k,\rho_k^{-1}(r)}}\frac{\partial\varphi_k}{\partial\nu}(y)d\sigma_{k,\rho_k^{-1}(r)}(y),$$
where $\nu$ denotes the outward unit normal to $\partial\Omega_{k,\rho_k^{-1}(r)}$ (note that $\frac{\partial\varphi_k}{\partial\nu}<0$ on $\partial\Omega_{k,\rho_k^{-1}(r)}$ by definition of $\Omega_{k,\rho_k^{-1}(r)}$). Green-Riemann formula and the choice of $u_k$ yield
\begin{equation}\label{eqsym1}\baa{rcl}
\displaystyle{\int_{\Sigma_{k,\rho_k^{-1}(r)}}}|\nabla\varphi_k(y)|d\sigma_{k,\rho_k^{-1}(r)}(y) & = & -\displaystyle{\int_{\Omega_{k,\rho_k^{-1}(r)}}}\Delta\varphi_k(y)dy\\
& = & -n\alpha_nr^{n-1}v_k(r)=n\alpha_nr^{n-1}|\nabla u_k(x)|.\eaa
\end{equation}\par
In the sequel, let
$$S_{s,s'}=\{z\in\R^n,\ s<|z|<s'\}$$
be the spherical shell of $\R^n$ between the radii $s$ and $s'$, with $0\le s<s'$. Let $t$ be any real number in $(0,\eta)$. A similar calculation as above gives
$$\int_{S_{r-t,r}}\Delta u_k(y)dy=\int_{\partial S_{r-t,r}}\frac{\partial u_k}{\partial\nu}(y)d\sigma(y)=n\alpha_n[r^{n-1}v_k(r)-(r-t)^{n-1}v_k(r-t)],$$
where $d\sigma$ and $\nu$ denote the superficial measure on $\partial S_{r-t,r}$ and the outward unit normal to $S_{r-t,r}$. By definition of $v_k$, one gets that
\begin{equation}\label{eqsym2}
\int_{S_{r-t,r}}\Delta u_k(y)dy=\int_{\Omega_{k,\rho_k^{-1}(r)}\backslash\Omega_{k,\rho_k^{-1}(r-t)}}\Delta\varphi_k(y)dy.
\end{equation}\par
For $a\in[0,M_k)$, call
$$j_k(a)=\int_{\Omega_{k,a}}|\nabla\varphi_k(y)|dy.$$
The Cauchy-Schwarz inequality gives
$$\baa{rcl}
[j_k(\rho_k^{-1}(r))-j_k(\rho_k^{-1}(r-t))]^2 & = & \left[\displaystyle{\int_{\Omega_{k,\rho_k^{-1}(r)}\backslash\Omega_{k,\rho_k^{-1}(r-t)}}}|\nabla\varphi_k(y)|dy\right]^2\\
& \le & \left|\Omega_{k,\rho_k^{-1}(r)}\backslash\Omega_{k,\rho_k^{-1}(r-t)}\right|\times\displaystyle{\int_{\Omega_{k,\rho_k^{-1}(r)}\backslash\Omega_{k,\rho_k^{-1}(r-t)}}}|\nabla\varphi_k(y)|^2dy.\eaa$$
Thus,
\begin{equation}\label{eqsym2bis}\baa{rl}
\displaystyle{\frac{\displaystyle{\int_{\Omega_{k,\rho_k^{-1}(r)}\backslash\Omega_{k,\rho_k^{-1}(r-t)}}}|\nabla\varphi_k(y)|dy}{\left|\Omega_{k,\rho_k^{-1}(r)}\backslash\Omega_{k,\rho_k^{-1}(r-t)}\right|}}\ \le\ A\ := & \displaystyle{\frac{\rho_k^{-1}(r-t)-\rho_k^{-1}(r)}{j_k(\rho_k^{-1}(r))-j_k(\rho_k^{-1}(r-t))}}\\
& \times\displaystyle{\frac{\displaystyle{\int_{\Omega_{k,\rho_k^{-1}(r)}\backslash\Omega_{k,\rho_k^{-1}(r-t)}}}|\nabla\varphi_k(y)|^2dy}{\rho_k^{-1}(r-t)-\rho_k^{-1}(r)}}.\eaa
\end{equation}
The first factor in the right-hand side of the above inequality converges to $1/i_k(\rho_k^{-1}(r))$ as $t\to 0^+$ (where $i_k$ is defined in (\ref{ghi})), from the co-area formula and the facts that $i_k$ is continuous in $Y_k\ni\rho_k^{-1}(r)$ and $\rho_k^{-1}$ is continuous in $[0,R]$. Similarly,
$$\frac{\displaystyle{\int_{\Omega_{k,\rho_k^{-1}(r)}\backslash\Omega_{k,\rho_k^{-1}(r-t)}}}|\nabla\varphi_k(y)|^2dy}{\rho_k^{-1}(r-t)-\rho_k^{-1}(r)}\ \mathop{\to}_{t\to 0^+}\ \int_{\Sigma_{k,\rho_k^{-1}(r)}}|\nabla\varphi_k(y)|d\sigma_{k,\rho_k^{-1}(r)}(y)=n\alpha_nr^{n-1}|\nabla u_k(x)|$$
from (\ref{eqsym1}). Therefore,
\begin{equation}\label{eqsym2ter}
A\ \mathop{\to}_{t\to 0^+}\ \frac{n\alpha_nr^{n-1}}{i_k(\rho_k^{-1}(r))}|\nabla u_k(x)|\ \le\ |\nabla u_k(x)|
\end{equation}
from the isoperimetric inequality applied to $\Sigma_{k,\rho_k^{-1}(r)}=\partial\Omega_{k,\rho_k^{-1}(r)}$ and $\partial B^n_r$ (namely, $n\alpha_nr^{n-1}\le i_k(\rho_k^{-1}(r))$). As a consequence,
\begin{equation}\label{eqsym3}
\mathop{\limsup}_{t\to 0^+}\frac{\displaystyle{\int_{\Omega_{k,\rho_k^{-1}(r)}\backslash\Omega_{k,\rho_k^{-1}(r-t)}}}|\nabla\varphi_k(y)|dy}{\left|\Omega_{k,\rho_k^{-1}(r)}\backslash\Omega_{k,\rho_k^{-1}(r-t)}\right|}\le|\nabla u_k(x)|.
\end{equation}\par
Let $\omega\ge 0$ be fixed and let $(\epsilon_l)_{l\in\N}$ be a sequence of positive real numbers, such that $\epsilon_l\to 0$ as $l\to+\infty$. It follows from (\ref{eqsym2}) and (\ref{eqsym3}) that there exists a sequence of positive numbers $t_l\in(0,\eta)$ such that $t_l\to 0$ and
\begin{equation}\label{eqsym4}\baa{l}
\displaystyle{\frac{\displaystyle{\int_{\Omega_{k,\rho_k^{-1}(r)}\backslash\Omega_{k,\rho_k^{-1}(r-t_l)}}}[\Delta\varphi_k(y)+\omega\ |\nabla\varphi_k(y)|]dy}{\left|\Omega_{k,\rho_k^{-1}(r)}\backslash\Omega_{k,\rho_k^{-1}(r-t_l)}\right|}}\\
\qquad\qquad\le\displaystyle{\frac{\displaystyle{\int_{S_{r-t_l,r}}}\Delta u_k(y)dy}{\left|\Omega_{k,\rho_k^{-1}(r)}\backslash\Omega_{k,\rho_k^{-1}(r-t_l)}\right|}}+\omega\ (1+\epsilon_l)|\nabla u_k(x)|.\eaa
\end{equation}
Since $u_k$ is radially symmetric and of class $C^2$ in $E_k$ (Lemma \ref{lem34}), and since
$$\left|\Omega_{k,\rho_k^{-1}(r)}\backslash\Omega_{k,\rho_k^{-1}(r-t_l)}\right|=|S_{r-t_l,r}|$$
for all $l\in\N$, the right-hand side of (\ref{eqsym4}) converges to $\Delta u_k(x)+\omega\ |\nabla u_k(x)|$ as $l\to+\infty$.\par
On the other hand, since the function $\varphi_k$ is of class $C^2(\overline{\Omega})$, (\ref{eqsym4}) shows that there exists a sequence of points $(y_l)_{l\in\N}$ of $\overline{\Omega}$ such that $\varphi_k(y_l)\in[\rho_k^{-1}(r),\rho_k^{-1}(r-t_l)]$ and
$$\mathop{\limsup}_{l\to+\infty}\ [\Delta\varphi_k(y_l)+\omega\ |\nabla\varphi_k(y_l)|]\le\Delta u_k(x)+\omega\ |\nabla u_k(x)|.$$
Up to extraction of some subsequence, one can assume that $y_l\to y\in\overline{\Omega}$, with $\varphi_k(y)=\rho_k^{-1}(r)$ (namely $y\in\Sigma_{k,\rho_k^{-1}(r)}$) and
$$\Delta\varphi_k(y)+\omega\ |\nabla\varphi_k(y)|\le\Delta u_k(x)+\omega\ |\nabla u_k(x)|.$$
That completes the proof of inequality (\ref{eq1pro35}).\par
For the proof of (\ref{eq2pro35}), let us first observe that the function $U_k$ is differentiable in $Y_k$, from Lemma \ref{lem32} and the fact that $u_k$ is of class $C^1(\overline{\Omega^*})$. Furthermore, since $u_k$ is radially symmetric and decreasing with respect to the variable $|x|$ and since $\rho_k$ is itself decreasing, it is enough to prove that
$$|\rho_k'(\rho_k^{-1}(|x|))|\times|\nabla u_k(x)|\ge 1$$
for all $x\in E_k$.\par
Fix $x\in E_k$ and use the same notations as above. It follows from (\ref{eqsym2bis}) that, for $t\in(0,\eta)$,
\begin{equation}\label{eqsym5}\baa{rcl}
1 & \le & A\times\displaystyle{\frac{\left|\Omega_{k,\rho_k^{-1}(r)}\backslash\Omega_{k,\rho_k^{-1}(r-t)}\right|}{\displaystyle{\int_{\Omega_{k,\rho_k^{-1}(r)}\backslash\Omega_{k,\rho_k^{-1}(r-t)}}}|\nabla\varphi_k(y)|dy}}\\
& = & A\times\displaystyle{\frac{\alpha_n[\rho_k(\rho_k^{-1}(r))^n-\rho_k(\rho_k^{-1}(r-t))^n]}{\rho_k^{-1}(r)-\rho_k^{-1}(r-t)}}\times\displaystyle{\frac{\rho_k^{-1}(r)-\rho_k^{-1}(r-t)}{\displaystyle{\int_{\Omega_{k,\rho_k^{-1}(r)}\backslash\Omega_{k,\rho_k^{-1}(r-t)}}}|\nabla\varphi_k(y)|dy}}.\eaa
\end{equation}
Since $\rho_k^{-1}$ is continuous in $[0,R]$ and $\rho_k$ is differentiable in $Y_k\ni\rho_k^{-1}(r)$, the second factor in the right-hand side of the above inequality converges to $n\alpha_nr^{n-1}\rho_k'(\rho_k^{-1}(r))=-n\alpha_nr^{n-1}|\rho_k'(\rho_k^{-1}(r))|$ as $t\to 0^+$. On the other hand, as already underlined earlier in the proof, the third factor converges to $-1/i_k(\rho_k^{-1}(r))$ as $t\to 0^+$. Together with (\ref{eqsym2ter}), the limit as $t\to 0^+$ in (\ref{eqsym5}) leads to
\begin{equation}\label{eqsym6}
1\le\frac{n\alpha_nr^{n-1}}{i_k(\rho_k^{-1}(r))}\times|\rho_k'(\rho_k^{-1}(r))|\times|\nabla u_k(x)|.
\end{equation}
The isoperimetric inequality yields $n\alpha_nr^{n-1}\le i_k(\rho_k^{-1}(r))$, whence
$$1\le|\rho_k'(\rho_k^{-1}(r))|\times|\nabla u_k(x)|$$
and the proof of Proposition \ref{pro35} is complete.\hfill\fin

\begin{lem}\label{lem37} For all $\epsilon>0$, there exists $k_0\in\N$ such that
\begin{equation}\label{inequk}
-\Delta u_k-(\tau+\epsilon)\ |\nabla u_k|\ \le\ [\underline{\lambda}(\Omega,\tau)+\epsilon]\ u_k\hbox{ in }E_k
\end{equation}
for all $k\ge k_0$.
\end{lem}

\noindent{\bf{Proof.}} Let us first recall that $\varphi$ is of class $C^{2,\alpha}(\overline{\Omega})$ and satisfies (\ref{eqminbis}). As already underlined at the beginning of this section, $|\nabla\varphi|$ and $\varphi$ are continuous nonnegative functions in $\overline{\Omega}$, which do not vanish simultaneously. There exists then $\gamma>0$ such that
$$|\nabla\varphi|+\varphi\ge\gamma>0\ \hbox{ in }\overline{\Omega}.$$\par
Fix $\epsilon>0$. Since
$$-\Delta\varphi_k=f_k\ \to\ f=\tau|\nabla\varphi|+\underline{\lambda}(\Omega,\tau)\varphi=-\Delta\varphi\ \hbox{ as }k\to+\infty$$
uniformly in $\overline{\Omega}$, standard elliptic estimates imply that $\varphi_k\to\varphi$ as $k\to+\infty$ in $W^{2,p}(\Omega)$ for any $1\le p<+\infty$, whence $\varphi_k\to\varphi$ in $C^{1,\beta}(\overline{\Omega})$ for all $\beta\in[0,1)$. As a consequence,
$$\baa{rcl}
\Delta\varphi_k+(\tau+\epsilon)|\nabla\varphi_k|+[\underline{\lambda}(\Omega,\tau)+\epsilon]\ \varphi_k & \displaystyle{\mathop{\to}_{k\to+\infty}} & \Delta\varphi+(\tau+\epsilon)|\nabla\varphi|+[\underline{\lambda}(\Omega,\tau)+\epsilon]\ \varphi\\
& & =\epsilon(|\nabla\varphi|+\varphi)\eaa$$
uniformly in $\overline{\Omega}$. But
$$\epsilon(|\nabla\varphi|+\varphi)\ge\epsilon\gamma>0 \hbox{ in }\overline{\Omega}$$
from the choice of $\gamma$. Therefore, there exists $k_0\in\N$ such that
\begin{equation}\label{k0}
\Delta\varphi_k+(\tau+\epsilon)|\nabla\varphi_k|+[\underline{\lambda}(\Omega,\tau)+\epsilon]\ \varphi_k\ge 0\ \hbox{ in }\overline{\Omega}
\end{equation}
for all $k\ge k_0$.\par
Let $k\in\N$ be such that $k\ge k_0$ and let $x$ be in $E_k$. From Proposition \ref{pro35}, there exists $y\in\overline{\Omega}$ such that $\varphi_k(y)=\rho_k^{-1}(|x|)$ and
$$\Delta u_k(x)+(\tau+\epsilon)|\nabla u_k(x)|\ge\Delta\varphi_k(y)+(\tau+\epsilon)|\nabla\varphi_k(y)|.$$
Thus,
\begin{equation}\label{taueps}
\Delta u_k(x)+(\tau+\epsilon)|\nabla u_k(x)|\ge-[\underline{\lambda}(\Omega,\tau)+\epsilon]\ \varphi_k(y)=-[\underline{\lambda}(\Omega,\tau)+\epsilon]\ \rho_k^{-1}(|x|)
\end{equation}
from (\ref{k0}). Inequality (\ref{inequk}) follows from Corollary \ref{cor36} and the fact that $\underline{\lambda}(\Omega,\tau)+\epsilon\ge 0$.\hfill\fin

\begin{lem}\label{lem38} For all $\epsilon>0$,
\begin{equation}\label{ineqlambda}
\underline{\lambda}(\Omega,\tau)+\epsilon\ge\lambda_1(\Omega^*,(\tau+\epsilon)e_r)=\underline{\lambda}(\Omega^*,\tau+\epsilon).
\end{equation}
\end{lem}

\noindent{\bf{Proof.}} Fix $\epsilon>0$ and let $k\in\N$ be large enough so that (\ref{inequk}) holds. Let $\psi$ be the solution (unique up to multiplication) of (\ref{eqball}) with coefficient $\tau+\epsilon$ instead of $\tau$. Namely, $\psi\in C^2(\overline{\Omega^*})$ solves
\begin{equation}\label{eqballbis}\left\{\baa{rcll}
-\Delta\psi-(\tau+\epsilon)|\nabla\psi|=-\Delta\psi+(\tau+\epsilon)e_r\cdot\nabla\psi & = & \lambda_1(\Omega^*,(\tau+\epsilon)e_r)\psi\\
& = & \underline{\lambda}(\Omega^*,\tau+\epsilon)\psi & \hbox{in }\Omega^*\\
\psi & > & 0 & \hbox{in }\Omega^*\\
\psi & = & 0 & \hbox{on }\partial\Omega^*.\eaa\right.
\end{equation}\par
Assume that the conclusion of Lemma \ref{lem38} does not hold, namely assume that
\begin{equation}\label{assume}
\underline{\lambda}(\Omega,\tau)+\epsilon<\underline{\lambda}(\Omega^*,\tau+\epsilon).
\end{equation}
We shall argue as in the proof of Lemma \ref{comparison} and compare $u_k$ with $\psi$. Let us first point out that, from Hopf lemma, $\frac{\partial\psi}{\partial\nu}$ is negative and continuous on the compact $\partial\Omega^*$, where $\nu$ is the
outward unit normal to $\partial\Omega^*$. Furthermore, $\psi$ is (at least) of class $C^1(\overline{\Omega^*})$ and is positive in $\Omega^*$. On the other hand, the function $u_k$ is (at least) of class $C^1(\overline{\Omega^*})$. Therefore, there exists $\eta_0>0$ such that
$$\forall\ \eta\ge\eta_0,\quad u_k\le\eta\psi\ \hbox{ in }\overline{\Omega^*}.$$\par
Call
$$\eta^*=\inf\ \{\eta>0,\ u_k\le\eta\psi\hbox{ in }\overline{\Omega^*}\}.$$
The real number $\eta^*$ is positive since $u_k>0$ in $\Omega^*$. There holds
$$u_k\le\eta^*\psi\hbox{ in }\overline{\Omega^*}.$$
Two cases may then occur~: either $u_k<\eta^*\psi$ in $\Omega^*$, or $\min_{\Omega^*}(\eta^*\psi-u_k)=0$. Let us first deal with\par
{\it Case 1}~: $u_k<\eta^*\psi$ in $\Omega^*$. Since $u_k$ is of class $C^1(\overline{\Omega^*}$), radially symmetric and decreasing with respect to the variable $|x|$, one has $e_r\cdot\nabla u_k=-|\nabla u_k|$ in $\overline{\Omega^*}$. It then follows from (\ref{inequk}) and (\ref{assume}) that
$$-\Delta u_k+(\tau+\epsilon)e_r\cdot\nabla u_k<\underline{\lambda}(\Omega^*,\tau+\epsilon)u_k\ \hbox{ in }E_k\cap\Omega^*$$
(remember that $u_k>0$ in $\Omega^*$).\par
Call $z=\eta^*\psi-u_k$. From the definition of $\psi$ in (\ref{eqballbis}), one gets that
\begin{equation}\label{ineqz}
-\Delta z+(\tau+\epsilon)e_r\cdot\nabla z>\underline{\lambda}(\Omega^*,\tau+\epsilon)z\ \hbox{ in }E_k\cap\Omega^*.
\end{equation}
But $z$ is positive in $\Omega^*$ and it is of class $C^2(E_k)$ and of class $C^1(\overline{\Omega^*})$. Furthermore, $z$ vanishes on $\partial\Omega^*=\{x\in\R^n,\ |x|=R\}$ and
$$E_k\cap\Omega^*\supset S_{\rho_k(a_{k,1}),R}=\{x,\ \rho_k(a_{k,1})<|x|<R\},$$
with $\rho_k(a_{k,1})<R$. Hopf lemma yields $\frac{\partial z}{\partial\nu}<0$ on $\partial\Omega^*$. As in the beginning of the proof of this lemma, there exists then $\epsilon_0>0$ such that $z\ge\epsilon'u_k$ in $\overline{\Omega^*}$ for all $\epsilon'\in[0,\epsilon_0]$. Hence, $u_k\le\eta\psi$ in $\overline{\Omega^*}$ for all $\eta\ge\eta^*/(1+\epsilon_0)$. That contradicts the minimality of $\eta^*$ and case 1 is ruled out.\par
{\it Case 2}~: the case where $\min_{\Omega^*}(\eta^*\psi-u_k)=\min_{\Omega^*}z=0$ is itself divided into two subcases~: either $\min_{\Omega^*\cap E_k}z=0$, or $z>0$ in $\Omega^*\cap E_k$ and there exists $y\in\Omega^*\backslash E_k$ such that $z(y)=0$.\par
In the first subcase, since $E_k\cap\Omega^*$ is open and $z$, which satisfies (\ref{ineqz}), is nonnegative and vanishes at some point $x\in E_k\cap\Omega^*$, the strong maximum principle implies that $z$ vanishes in the connected component of $E_k\cap\Omega^*$ containing $x$. That is impossible because of the strict inequality in (\ref{ineqz}).\par
Therefore, only the second subcase could occur. In that subcase, owing to the definition of $E_k$, $|y|=\rho_k(a_{k,i})<R$ for some $i\in\{1,\cdots,m_k\}$ and there is $r_0>0$ such that
$$S=S_{|y|,|y|+r_0}=\{y',\ |y|<|y'|<|y|+r_0\}\subset E_k\cap\Omega^*.$$
The function $z$ is of class $C^2(S)\cap C^1(\overline{S})$, it is positive in $S$ and vanishes at $y\in\partial S$. Furthermore, $z$ satisfies (\ref{ineqz}) (at least) in $S$. Hopf lemma implies that $e_r\cdot\nabla z(y)>0$, where $e_r=y/|y|$ if $y\neq 0$, and $e_r$ may be any unit vector if $y=0$. But $z$ is of class $C^1(\overline{\Omega^*})$ and it has a minimum at $y\in\Omega^*$, whence $\nabla z(y)=0$. This gives a contradiction.\par
Case 2 is then ruled out too and the proof of Lemma \ref{lem38} is complete.\hfill\fin\break

\noindent{\bf{Proof of the inequality $\lambda_1(\Omega,v)\ge\lambda_1(\Omega^*,\tau e_r)$.}} Let us now complete the proof of the first part of Theorem \ref{th1}. To do so, remember that the function $\omega\mapsto\underline{\lambda}(\Omega^*,\omega)$ is continuous in $\R_+$ (see Section \ref{secfurther}, Lemma \ref{continuity}). Therefore, passing to the limit as $\epsilon\to 0$ in (\ref{ineqlambda}) yields
$$\underline{\lambda}(\Omega,\tau)\ge\underline{\lambda}(\Omega^*,\tau).$$
On the other hand, from Theorem \ref{th3}, $\underline{\lambda}(\Omega^*,\tau)=\lambda_1(\Omega^*,\tau e_r)$. As a conclusion,
$$\lambda_1(\Omega,v)\ge\underline{\lambda}(\Omega,\tau)\ge\underline{\lambda}(\Omega^*,\tau)=\lambda_1(\Omega^*,\tau e_r)$$
for all $v\in L^{\infty}(\Omega)$ with $\|v\|_{L^{\infty}(\Omega)}\le\tau$. That completes the proof of formula (\ref{main}) in Theorem~\ref{th1}.\hfill\fin


\subsection{Characterization of the case of equality in (\ref{main})}\label{sec32}

We use in this section the same notations as in Section \ref{sec31}. Assume first that $\Omega$ is a ball, say with the origin as center (up to translation). In other words, assume that $\Omega=\Omega^*$. It follows from Theorems \ref{th2} and \ref{th3} that, for all $v\in L^{\infty}(\Omega^*)$ with $\|v\|_{L^{\infty}(\Omega^*)}\le\tau$, the equality $\lambda_1(\Omega^*,v)=\lambda_1(\Omega^*,\tau e_r)$ holds only when $v=\tau e_r$.\par
Consider now the case where $\Omega$ is not a ball and call $R=(|\Omega|/\alpha_n)^{1/n}$ the radius of $\Omega^*$. We shall prove that
$$\lambda_1(\Omega,v)>\lambda_1(\Omega^*,\tau e_r)=\underline{\lambda}(\Omega^*,\tau)$$
for all $v\in L^{\infty}(\Omega)$ with $\|v\|_{L^{\infty}(\Omega)}\le\tau$. Owing to the definition of
$\underline{\lambda}(\Omega,\tau)$, it is enough to prove that
\begin{equation}\label{main2}
\underline{\lambda}(\Omega,\tau)>\underline{\lambda}(\Omega^*,\tau).
\end{equation}\par
The proof is divided into several lemmata.\par
First, the isoperimetric inequality yields the existence of $\beta>0$ such that
\be\label{areabeta}
\hbox{area}(\partial\Omega)=\int_{\partial\Omega}d\sigma_{\partial\Omega}(y)\ge(1+\beta)n\alpha_nR^{n-1},
\ee
where the left-hand side is the $(n-1)$-dimensional area of $\partial\Omega$.\par
Call $d(x,{\mathcal{A}})$ the Euclidean distance of a point $x\in\R^n$ to a set $\mathcal{A}\subset\R^n$. For all $\gamma>0$, define
$$U_{\gamma}=\{x\in\overline{\Omega},\ d(x,\partial\Omega)\le\gamma\}$$
and for all $y\in\partial\Omega$, call $\nu(y)$ the outward unit normal to $\Omega$. Since $\partial\Omega$ is of class $C^2$, there exists $\gamma_1>0$ such that the segments $[y,y-\gamma_1\nu(y)]$ are included in
$\overline{\Omega}$ and pairwise disjoint when $y$ describes $\partial\Omega$ (thus, the ``segments'' $(y,y-\gamma_1\nu(y)]$ describe the set $\{x\in\Omega$, $d(x,\partial\Omega)\le\gamma_1\}$ as $y$ describes $\partial\Omega$).

\begin{lem}\label{lem39} Let $\varphi\in C^{2,\alpha}(\overline{\Omega})$ solve (\ref{eqminbis}) with
$\|\varphi\|_{L^{\infty}(\Omega)}=1$ and call
$$m=\mathop{\min}_{y\in\partial\Omega}\ |\nabla\varphi(y)|.$$
Then $m>0$ and there exists $\gamma_2\in(0,\gamma_1]$ such that, for all $\gamma\in(0,\gamma_2]$, $|\nabla\varphi|\neq 0$ in
$U_{\gamma}$, $\nabla\varphi(y-r\nu(y))\cdot\nu(y)<0$ for all $y\in\partial\Omega$ and $r\in[0,\gamma]$, and $\varphi\ge\gamma
m/2$ in $\overline{\Omega\backslash U_{\gamma}}$.
\end{lem}

\noindent{\bf{Proof.}} Let us first observe that $m>0$ since $\varphi$ is (at least) of class $C^1(\overline{\Omega})$ and
$\frac{\partial\varphi}{\partial\nu}(y)=\nabla\varphi(y)\cdot\nu(y)<0$ for all $y\in\partial\Omega$ (from Hopf lemma).\par
Assume that the conclusion of the lemma does not hold. Then there exists a sequence of positive numbers $(\gamma^l)_{l\in\N}\to
0$ such that one of the three following cases occur~: 1) either for each $l\in\N$, there is a point $x_l\in U_{\gamma^l}$ such
that $\nabla\varphi(x_l)=0$, 2) or for each $l\in\N$, there are a point $y_l\in\partial\Omega$ and a number $r_l\in[0,\gamma^l]$
such that $\nabla\varphi(y_l-r_l\nu(y_l))\cdot\nu(y_l)\ge 0$, 3) or for each $l\in\N$, there is a point
$x_l\in\overline{\Omega\backslash U_{\gamma^l}}$ such that $\varphi(x_l)<\gamma^lm/2$.\par
In the first case, after passing to the limit up to extraction of some subsequence, there would exist a point
$x\in\partial\Omega$ such that $\nabla\varphi(x)=0$. This is impossible. Similarly, in the second case, there would exist a
point $y\in\partial\Omega$ such that $\nabla\varphi(y)\cdot\nu(y)\ge 0$, which is still impossible.\par
Assume that the third case occurs. Let $y_l\in\partial\Omega$ be such that
$$d_l:=|x_l-y_l|=d(x_l,\partial\Omega)\ge\gamma^l>0.$$
Up to extraction of some subsequence, one has $x_l\to x\in\overline{\Omega}$ and $\varphi(x)\le 0$ by passing to the limit as
$l\to+\infty$ in the inequality $\varphi(x_l)<\gamma^lm/2$. Since $\varphi>0$ in $\Omega$, it follows that $x\in\partial\Omega$,
whence $|x_l-y_l|\to 0$ and $y_l\to x$ as $l\to+\infty$. On the one hand, the mean value theorem implies that
$$\frac{\varphi(x_l)-\varphi(y_l)}{|x_l-y_l|}=\nabla\varphi(z_l)\cdot\frac{x_l-y_l}{|x_l-y_l|}\to-\nabla\varphi(x)\cdot\nu(x)=|\nabla\varphi(x)|\
\hbox{ as }l\to+\infty,$$
where $z_l$ is a point lying on the segment between $x_l$ and $y_l$ (whence, $z_l\to x$ as $l\to+\infty$). On the other hand,
since $\varphi=0$ on $\partial\Omega$,
$$\frac{\varphi(x_l)-\varphi(y_l)}{|x_l-y_l|}=\frac{\varphi(x_l)}{d_l}<\frac{\gamma^lm}{2\gamma^l}=\frac{m}{2}.$$
Hence, $|\nabla\varphi(x)|\le m/2$ at the limit, which contradicts the definition and the positivity of $m$.\hfill\fin\break

In the sequel, we use the same functions $\varphi_k$ as in Section \ref{sec31}, together with the same sets $Z_k$, $Y_k$,
$\Omega_{k,a}$, $\Sigma_{k,a}$ and functions $\rho_k$, $u_k$, etc.

\begin{lem}\label{lem310} There exist $k_1\in\N$ and $a_0>0$ such that $[0,a_0]\subset Y_k$ for all $k\ge k_1$, and 
$$i_k(a)=\int_{\Sigma_{k,a}}d\sigma_{k,a}(y)=\hbox{area}(\Sigma_{k,a})\ge\left(1+\frac{\beta}{2}\right)n\alpha_nR^{n-1}$$
for all $k\ge k_1$ and $a\in[0,a_0]$, where $\beta>0$ is given in (\ref{areabeta}).
\end{lem}

\noindent{\bf{Proof.}} Let $\gamma_2>0$ be as in Lemma \ref{lem39}. By compactness of $U_{\gamma_2}$ and $\partial\Omega$, and since $\varphi$ is of class
$C^1(\overline{\Omega})$, there is $\delta>0$ such that $|\nabla\varphi|\ge\delta$ in $U_{\gamma_2}$ and
$\nabla\varphi(y-r\nu(y))\cdot\nu(y)\le-\delta$ for all $y\in\partial\Omega$ and $r\in[0,\gamma_2]$.\par
Since $\varphi_k\to\varphi$ in $C^1(\overline{\Omega})$ as $k\to+\infty$ (the convergence actually holds in
$C^{1,\alpha'}(\overline{\Omega})$ for all $0\le\alpha'<1$), there exists $k_1\in\N$ such that
$$\forall\ k\ge k_1,\quad|\nabla\varphi_k|\ge\delta/2>0\hbox{ in }U_{\gamma_2},\quad\varphi_k\ge\gamma_2 m/4>0\hbox{ in }\overline{\Omega\backslash U_{\gamma_2}}$$
and
\be\label{phikdelta}
\forall\ k\ge k_1,\ \forall\ y\in\partial\Omega,\ \forall\ r\in[0,\gamma_2],\quad\nabla\varphi_k(y-r\nu(y))\cdot\nu(y)\le-\delta/2<0.
\ee\par
Let $k\in\N$ be fixed such that $k\ge k_1$. It especially follows that, for all $a\in[0,\gamma_2 m/8]$,
\begin{equation}\label{sigmaka}
\Sigma_{k,a}=\{x\in\overline{\Omega},\ \varphi_k(x)=a\}\subset U_{\gamma_2},
\end{equation}
whence $|\nabla\varphi_k|\neq 0$ everywhere on the $C^2$ hypersurface $\Sigma_{k,a}$. Thus,
$$[0,\gamma_2 m/8]\subset Y_k.$$
Furthermore, for all $y\in\partial\Omega$, the segment $[y,y-\gamma_2\nu(y)]$ is included in $\overline{\Omega}$
and there exists $\theta\in[0,1]$ such that
$$\varphi_k(y-\gamma_2\nu(y))=\underbrace{\varphi_k(y)}_{=0}-\gamma_2\nu(y)\cdot\nabla\varphi_k(y-\theta\gamma_2\nu(y))\ge\frac{\gamma_2\delta}{2}~;$$
more precisely, the function $\kappa\ :\ [0,\gamma_2]\to\R$, $s\mapsto\varphi_k(y-s\nu(y))$ is differentiable and
$\kappa'(s)\ge\delta/2$ for all $s\in[0,\gamma_2]$.\par
Call
$$a_1=\min(\gamma_2 m/8,\gamma_2\delta/4)>0.$$
It follows from the above calculation that, for all $k\ge k_1$, $a\in[0,a_1]$ and $y\in\partial\Omega$, there exists a unique
point $\phi_{k,a}(y)\in[y,y-\gamma_2\nu(y)]\cap\Sigma_{k,a}$. Moreover, for such a choice of $k$ and $a$, the map $\phi_{k,a}$ is one-to-one since the segments $[y,y-\gamma_2\nu(y)]$ are pairwise
disjoint when $y$ describes $\partial\Omega$ (because $\gamma_2\in(0,\gamma_1]$). Lastly,
\be\label{onto}
\Sigma_{k,a}=\{\phi_{k,a}(y),\ y\in\partial\Omega\}
\ee
from (\ref{sigmaka}).\par
Let us now prove that the area of $\Sigma_{k,a}$ is close to that of $\partial\Omega$ for $k$ large enough and $a\ge 0$ small enough. To do so, let us first represent $\partial\Omega$ by a finite number of $C^{2,\alpha}$ maps $y^1,\ldots,y^p$ (for some $p\in\N^*$) defined in $\overline{\mathcal{B}}$, and for which
$$\partial_1y^j(x')\times\cdots\times\partial_{n-1}y^j(x')\neq 0\ \hbox{ for all }1\le j\le p\hbox{ and }x'\in\overline{\mathcal{B}_{\rho}}.$$
Here, ${\mathcal{B}}=\{x'=(x_1,\ldots,x_{n-1}),\ |x'|<1\}$ and $\partial_iy^j(x')=(\partial_{x_i}y^j_1(x'),\ldots,\partial_{x_i}y^j_n(x'))$ for $1\le i\le n-1$, where $y^j(x')=(y^j_1(x'),\ldots,y^j_n(x'))\in\R^n$. The maps $y^j$ are chosen so that
$$\partial\Omega=\{y^j(x'),\ x'\in\overline{\mathcal{B}},\ 1\le j\le p\}.$$\par
Fix $k\ge k_1$ and $a\in[0,a_1]$. For each $1\le j\le p$, there exists then a map $t^j_{k,a}:\overline{\mathcal{B}}\to[0,\gamma_2]$ such that
\be\label{tjka}
\varphi_k(y^j(x')-t^j_{k,a}(x')\nu(y^j(x')))=a
\ee
for all $x'\in\overline{\mathcal{B}}$, and $\Sigma_{k,a}=\{y^j(x')-t^j_{k,a}(x')\nu(y^j(x')),\ x'\in\overline{\mathcal{B}},\ 1\le j\le p\}$. From the arguments above, each real number $t^j_{k,a}(x')$ is then uniquely determined, and $t^j_{k,0}(x')=0$.\par
Since the functions $\varphi_k$ (say, for all $k\ge k_1$), $y^j$ and $\nu\circ y^j$ (for all $1\le j\le p$) are (at least) of class $C^1$ (respectively in $\overline{\Omega}$, $\overline{\mathcal{B}}$ and $\overline{\mathcal{B}}$), it follows from implicit function theorem and (\ref{phikdelta}) that the functions $t_{k,a}^j$ (for all $k\ge k_1$, $a\in[0,a_1]$, $1\le j\le p$) and
$$h^j_{k,x'}:[0,a_1]\ni a\mapsto t_{k,a}^j(x')$$
(for all $k\ge k_1$, $1\le j\le p$, $x'\in\overline{\mathcal{B}}$) are of class $C^1$ (respectively in $\overline{\mathcal{B}}$ and $[0,a_1]$). From the chain rule applied to (\ref{tjka}), it is straightforward to check that, for all $k\ge k_1$, $a\in[0,a_1]$, $1\le j\le p$ and $x'\in\overline{\mathcal{B}}$,
$$(h^j_{k,x'})'(a)=\frac{-1}{\nu(y^j(x'))\ \cdot\ \nabla\varphi_k(y^j(x')-t^j_{k,a}(x')\nu(y^j(x')))}\in(0,2\delta^{-1}]\ \hbox{ from (\ref{phikdelta})},$$
whence
\be\label{tjka2}
0\le t^j_{k,a}(x')=h^j_{k,x'}(a)\le 2\delta^{-1}a
\ee
because $h^j_{k,x'}(0)=t^j_{k,0}(x')=0$.\par
Similarly,
\be\label{chain}
\partial_{x_i}t^j_{k,a}(x')=\frac{[\partial_iy^j(x')-t^j_{k,a}(x')\partial_i(\nu\circ y^j)(x')]\ \cdot\ \nabla\varphi_k(y^j(x')-t^j_{k,a}(x')\nu(y^j(x')))}{\nu(y^j(x'))\ \cdot\ \nabla\varphi_k(y^j(x')-t^j_{k,a}(x')\nu(y^j(x')))}
\ee
for all $k\ge k_1$, $a\in[0,a_1]$, $1\le j\le p$, $x'\in\overline{\mathcal{B}}$ and $1\le i\le n-1$. For all $k$, $1\le j\le p$ and $x'\in\overline{\mathcal{B}}$, one has $\varphi_k(y^j(x'))=0$, whence $\partial_iy^j(x')\cdot\nabla\varphi_k(y^j(x'))=0$ (for all $1\le i\le n-1$). On the other hand, the functions $\varphi_k$ converge (at least) in $C^{1,1/2}(\overline{\Omega})$ to $\varphi$ as $k\to+\infty$. As a consequence, 
$$|\partial_iy^j(x')\cdot\nabla\varphi_k(y^j(x')-t^j_{k,a}(x')\nu(y^j(x')))|\le C_1\sqrt{t^j_{k,a}(x')}$$
for all $k\ge k_1$, $a\in[0,a_1]$, $1\le j\le p$, $x'\in\overline{\mathcal{B}}$ and $1\le i\le n-1$, and for some constant $C_1$ defined by
$$C_1=\max_{1\le i'\le n-1,\ 1\le j'\le p,\ \xi\in\overline{\mathcal{B}}}|\partial_{i'}y^{j'}(\xi)|\ \times\ \sup_{k'\in\N,\ z\neq z'\in\Omega}\frac{|\nabla\varphi_{k'}(z)-\nabla\varphi_{k'}(z')|}{\sqrt{|z-z'|}}\ <\ +\infty.$$
Call now
$$C_2=\max_{1\le j'\le p,\ \xi\in\overline{\mathcal{B}},\ 1\le i'\le n-1}|\partial_{i'}(\nu\circ y^{j'})(\xi))|\ \times\ \sup_{k'\in\N,\ z\in\overline{\Omega}}|\nabla\varphi_{k'}(z)|\ <\ +\infty.$$
Together with (\ref{chain}) and (\ref{phikdelta}), the above arguments imply that
$$\baa{rcl}
|\partial_{x_i}t^j_{k,a}(x')| & \le & 2\delta^{-1}(C_1\sqrt{t^j_{k,a}(x')}+C_2t^j_{k,a}(x'))\\
& \le & 2\delta^{-1}(C_1\sqrt{2\delta^{-1}a}+2C_2\delta^{-1}a)\ \hbox{ from (\ref{tjka2})}\eaa$$
for all $k\ge k_1$, $a\in[0,a_1]$, $1\le j\le p$, $x'\in\overline{\mathcal{B}}$ and $1\le i\le n-1$.\par
It follows that
\be\label{suptjka}
\sup_{k\ge k_1,\ 1\le j\le p,\ x'\in\overline{\mathcal{B}},\ 1\le i\le n-1}|\partial_{x_i}(t^j_{k,a}\ \nu\circ y^j)(x')|\ \to\ 0\ \hbox{ as }a\to 0,\ 0\le a\le a_1.
\ee\par
On the other hand, there are some open sets $U^1,\ldots,U^p$ of ${\mathcal{B}}$ such that
$$\hbox{area}(\partial\Omega)=\sum_{j=1}^p\int_{U^j}|\partial_1y^j(x')\times\cdots\times\partial_{n-1}y^j(x')|\ dx',$$
where the sets $\{y^j(x'),\ x'\in U^j\}$ for $j=1,\ldots,p$ are pairwise disjoint and, for any $\epsilon>0$, there are some measurable sets $V^1,\ldots,V^p\subset\overline{\mathcal{B}}$ such that $\{y^j(x'),\ x'\in V^j,\ 1\le j\le p\}=\partial\Omega$, $V^j\supset U^j$ and $\displaystyle{\int_{\mathcal{B}}}{\bf{1}}_{V^j\backslash U^j}(x')dx'\le\epsilon$ for all $1\le j\le p$. Since all functions $y^j$ and $t^j_{k,a}\nu\circ y^j$ (for all $k\ge k_1$, $a\in[0,a_1]$, $1\le j\le p$) are of class $C^1$ in $\overline{\mathcal{B}}$, since each function $\phi_{k,a}$ is one-to-one and since (\ref{onto}) holds, it follows that
$$\hbox{area}(\Sigma_{k,a})=\sum_{j=1}^p\int_{U^j}|\partial_1(y^j-t^j_{k,a}\nu\circ y^j)(x')\times\cdots\times\partial_{n-1}(y^j-t^j_{k,a}\nu\circ y^j)(x')|\ dx'$$
for all $k\ge k_1$ and $a\in[0,a_1]$.\par
One concludes from (\ref{suptjka}) that
$$\sup_{k\ge k_1}\ |\hbox{area}(\Sigma_{k,a})-\hbox{area}(\partial\Omega)|\to 0\hbox{ as }a\to 0\hbox{ with }0\le a\le a_1.$$
Because $\beta$ in (\ref{areabeta}) is positive, there exists then $a_0\in(0,a_1]$ such that
$$i_k(a)=\hbox{area}(\Sigma_{k,a})\ge\left(1+\frac{\beta}{2}\right)n\alpha_nR^{n-1}$$
for all $k\ge k_1$ and $a\in[0,a_0]$.\par
That completes the proof of Lemma \ref{lem310}.\hfill\fin

\begin{lem}\label{lem311} With the notations of Lemma \ref{lem310}, one has
$$u_k(x)\ge\left(1+\frac{\beta}{2}\right)\rho_k^{-1}(|x|)$$
for all $k\ge k_1$ and $x\in\overline{\Omega^*}$ such that $\rho_k(a_0)\le|x|\le R$.
\end{lem}

\noindent{\bf{Proof.}} Fix $k\ge k_1$. From Lemma \ref{lem310}, one has $\overline{S_{\rho_k(a_0),R}}\subset E_k$. Fix any
$x\in\overline{\Omega^*}$ such that $r=|x|\in[\rho_k(a_0),R]$ (notice that $0<\rho_k(a_0)<R$ and $\rho_k^{-1}(r)\in[0,a_0]$).
The calculations of the proof of Proposition \ref{pro35}, and especially inequality (\ref{eqsym6}), imply that
$$1\le\frac{n\alpha_nr^{n-1}}{i_k(\rho_k^{-1}(r))}\times|\rho_k'(\rho_k^{-1}(r))|\times|\nabla u_k(x)|.$$
But 
$$\baa{rcl}
i_k(\rho_k^{-1}(r))=\displaystyle{\int_{\Sigma_{k,\rho_k^{-1}(r)}}}d\sigma_{k,\rho_k^{-1}(r)}(y)=\hbox{area}(\Sigma_{k,\rho_k^{-1}(r)})
& \ge & \left(1+\displaystyle{\frac{\beta}{2}}\right)n\alpha_nR^{n-1}\\
& \ge & \left(1+\displaystyle{\frac{\beta}{2}}\right)n\alpha_nr^{n-1}\eaa$$
from Lemma \ref{lem310}. Thus,
$$1+\frac{\beta}{2}\le|\rho_k'(\rho_k^{-1}(r))|\times|\nabla u_k(x)|.$$\par
The conclusion of Lemma \ref{lem311} follows from the above inequality, as in the proof of Corollary \ref{cor36}.\hfill\fin

\begin{lem}\label{lem312} There exist $k_2\ge k_1$ and $\eta>0$ such that
$$u_k(x)\ge(1+\eta)\ \rho_k^{-1}(|x|)$$
for all $k\ge k_2$ and $x\in\overline{\Omega^*}$.
\end{lem}

\noindent{\bf{Proof.}} Let $e$ be any unit vector in $\R^n$ and choose $k\ge k_1$. Let $\tilde{u}_k$ be the function defined in
$[0,R]$ by $\tilde{u}_k(r)=u_k(re)$ for all $r\in[0,R]$. This function is differentiable and decreasing in $[0,R]$. Furthermore,
Proposition \ref{pro35} and the fact that $\rho_k^{-1}$ is decreasing in $[0,R]$ imply that
$$-\tilde{u}_k'(r)\ge-\frac{d}{dr}(\rho_k^{-1}(r))$$
for all $r\in\rho_k(Y_k)=(0,\rho_k(a_{k,m-1}))\cup\cdots\cup(\rho_k(a_{k,2}),\rho_k(a_{k,1}))\cup(\rho_k(a_{k,1}),R]$. Finite
Increment Theorem yields especially, as in the proof of Corollary \ref{cor36},
$$u_k(re)-u_k(\rho_k(a_0)e)\ge\rho_k^{-1}(r)-a_0$$
for all $r\in[0,\rho_k(a_0)]$. Fix such a $r$ in $[0,\rho_k(a_0)]$ (whence $\rho_k^{-1}(r)\in[a_0,M_k]\subset(0,M_k]$). One gets
that
$$\frac{u_k(re)}{\rho_k^{-1}(r)}\ge 1+\frac{u_k(\rho_k(a_0)e)-a_0}{\rho_k^{-1}(r)}\ge 1+\frac{\beta a_0}{2\rho_k^{-1}(r)}$$
from Lemma \ref{lem311}.\par
But $\rho_k^{-1}(r)\le M_k=\max_{\overline{\Omega}}\varphi_k\to\max_{\overline{\Omega}}\varphi=1$ as $k\to+\infty$. Hence, there
exists $k_2\ge k_1$ such that
$$u_k(re)\ge\left(1+\frac{\beta a_0}{4}\right)\rho_k^{-1}(r)$$
for all $k\ge k_2$ and $r\in[0,\rho_k(a_0)]$. As in the proof of Corollary \ref{cor36}, the conclusion of Lemma \ref{lem312}
follows from the above inequality and from Lemma \ref{lem311}, with the choice
$$\eta=\min(\beta/2,\beta a_0/4)>0$$
for instance.\hfill\fin\break

\noindent{\bf{Conclusion.}} Fix any $\epsilon>0$. From Lemma \ref{lem312} and from (\ref{taueps}), there exists $k_3\ge k_2$
such that
$$-\Delta u_k(x)-(\tau+\epsilon)\ |\nabla u_k(x)|\ \le\ [\underline{\lambda}(\Omega,\tau)+\epsilon]\ \rho_k^{-1}(|x|)\ \le\
\frac{\underline{\lambda}(\Omega,\tau)+\epsilon}{1+\eta}\ u_k(x)$$
for all $x\in E_k$ and $k\ge k_3$. As in the proof of Lemma \ref{lem38}, one gets that
$$\frac{\underline{\lambda}(\Omega,\tau)+\epsilon}{1+\eta}\ge\underline{\lambda}(\Omega^*,\tau+\epsilon).$$
Passing to the limit $\epsilon\to 0^+$ in the above inequality yields
$$\underline{\lambda}(\Omega,\tau)\ge(1+\eta)\
\underline{\lambda}(\Omega^*,\tau)>\underline{\lambda}(\Omega^*,\tau)=\lambda_1(\Omega^*,\tau e_r)\ (>0).$$
That completes the proof of Theorem \ref{th1}.\hfill\fin


\SE{Appendix}\label{secapp}

                                                                           
\subsection{Behavior of $F_n(m,\tau)$ for large $\tau$}\label{secFn}

This section is devoted to the proof of the results mentioned in Remark \ref{rem2}.

First, to prove (\ref{F1}), fix $m>0$ and $\tau\geq 0$, set $\Omega=(-R,R)$ with $2R=m$, and
denote
$$\lambda=\lambda_1(\Omega, \tau e_r)$$
and $\varphi=\varphi_{\Omega,\tau e_r}$. Theorem \ref{th3} ensures that $\varphi$
is an even function, decreasing in $[0,R]$ and Theorem \ref{th2} yields
\[
-\varphi^{\prime\prime}(r)+\tau\varphi^{\prime}(r)=\lambda\varphi(r)\mbox{ for all }0\leq r\le R,
\]
with $\varphi(R)=0$, $\varphi>0$ in $(-R,R)$ and $\varphi^{\prime}(0)=0$. For all $s\in [0,\tau R]$, define $\psi(s)=\varphi(s/\tau)$, so that $\psi$
satisfies the equation
\[
-\psi^{\prime\prime}(s)+\psi^{\prime}(s)=\frac{\lambda}{\tau^2}\psi(s)\mbox{ for all }0\leq s\le\tau R,
\]
with $\psi(\tau R)=0$ and $\psi^{\prime}(0)=0$. Notice that $\lambda$ depends on $\tau$, but since, for all $\tau\geq 0$,
$0<\lambda\leq \lambda_1((-R,R),0)$, there exists $\tau_0>0$ such that $\tau^2\geq 4\lambda$ for all $\tau\geq \tau_0$, and we
will assume that $\tau\geq \tau_0$ in the sequel. The function $\psi$ can be computed explicitly: there exist $A,B\in \R$ such
that, for all $0\leq s\leq \tau R$,
\[
\psi(s)=Ae^{\mu_+r}+Be^{\mu_-r},
\]
where $\displaystyle \mu_{\pm}=(1\pm \sqrt{1-4\lambda/\tau^2})/2$. Using the boundary values of $\psi$ and
$\psi^{\prime}$, one obtains after straightforward computations:
\[
\lambda=\frac{\tau^2}4 \left(1+\sqrt{1-\frac{4\lambda}{\tau^2}}\right)^2 e^{-\sqrt{1-\frac{4\lambda}{\tau^2}}\tau R}.
\]
Since $\lambda$ remains bounded when $\tau\rightarrow +\infty$, it is then straightforward to check that $\lambda\sim \tau^2e^{-\tau R}$ when $\tau\rightarrow +\infty$, and that (\ref{F1}) follows.\hfill\break\par

We now turn to the proof of assertion (\ref{Fntau}). Let $n\geq 2$, $m>0$, $\tau\geq 0$ and $\Omega=B^n_R$ be such that $\left\vert \Omega\right\vert=m$, so that one has $R=(m/\alpha_n)^{1/n}$ and $F_n(m,\tau)=\lambda_1(\Omega,\tau e_r)$. We first claim that
\[
F_n(m,\tau)> F_1(2R,\tau).
\]
Indeed, write
$$\lambda=\lambda_1(\Omega,\tau e_r)\ \hbox{ and }\varphi_n=\varphi_{\Omega,\tau e_r}.$$
Similarly,
$F_1(2R,\tau)=\lambda_1((-R,R),\tau e_r)$, and we denote $\mu=\lambda_1((-R,R),\tau e_r)$ and $\varphi_1=\varphi_{(-R,R),\tau
e_r}$. As before, define $\psi_n(y)=\varphi_n(y/\tau)$ for all $y\in \tau\overline{\Omega}=\overline{B^n_{\tau R}}$ and $\psi_1(r)=\varphi_1(r/\tau)$ for all
$r\in[-\tau R,\tau R]$. Finally, since $\psi_n$ is radial, let $u_n:\left[0,\tau R\right]\rightarrow \R$ such that
$\psi_n(y)=u_n(\left\vert y\right\vert)$ for all $y\in \tau\overline{\Omega}=\overline{B^n_{\tau R}}$. One has
\begin{equation} \label{onedim}
\left\{
\begin{array}{ll}
\displaystyle -u_n^{\prime\prime}(r)-\frac{n-1}ru_n^{\prime}(r)+u_n^{\prime}(r)=\frac{\lambda}{\tau^2}u_n(r)& \mbox{ in }(0,\tau
R],\\
\\
\displaystyle -\psi_1^{\prime\prime}(r)+\psi_1^{\prime}(r)=\frac{\mu}{\tau^2}\psi_1(r) &\mbox{ in }[0,\tau R],
\end{array}
\right.
\end{equation}
with $u_n^{\prime}(0)=u_n(\tau R)=0$, $\psi^{\prime}_1(0)=\psi_1(\tau R)=0$. 

Assume that $\lambda\leq \mu$. Since $u_n^{\prime}<0$ in $(0,\tau R]$ and $u_n\ge 0$, one obtains
\be\label{onedimbis}
\left\{
\begin{array}{ll}
\displaystyle -u_n^{\prime\prime}(r)+u_n^{\prime}(r)\leq \frac{\mu}{\tau^2}u_n(r) &\mbox{ in }[0,\tau R],\\
\\
\displaystyle -\psi_1^{\prime\prime}(r)+\psi_1^{\prime}(r)=\frac{\mu}{\tau^2}\psi_1(r) &\mbox{ in }[0,\tau R].
\end{array}
\right.
\ee
Since $\psi_1^{\prime}(\tau R)<0$ by Hopf lemma, while $\psi_1(r)>0$ in $[0,\tau R)$, $u_n(r)>0$ in $[0,\tau R)$ and the functions
$u_n$ and $\psi_1$ belong (at least) to $C^1([0,\tau R])$, there exists then $\gamma>0$ such that
$\gamma\psi_1(r)>u_n(r)$ for all $0\le r<\tau R$. Define $\gamma^{\ast}$ as the infimum of all the $\gamma>0$ such that
$\gamma\psi_1>u_n$ in $[0,\tau R)$, observe that $\gamma^{\ast}>0$ and define $z=\gamma^{\ast}\psi_1-u_n$ which is non-negative
in $[0,\tau R]$ and satisfies
\begin{equation} \label{ineqdimone}
-z^{\prime\prime}(r)+z^{\prime}(r)-\frac{\mu}{\tau^2}z(r)\ge 0
\end{equation}
for all $0\le r\le\tau R$ and $z(\tau R)=0$. 

Assume that there exists $0<r<\tau R$ such that $z(r)=0$. The strong maximum principle shows that $z$ is identically zero in
$[0,\tau R]$, which means that $\gamma^{\ast}\psi_1=u_n$ in $[0,\tau R]$, and even that $\psi_1=u_n$ because
$\psi_1(0)=u_n(0)=1$. But this is impossible according to (\ref{onedim}) and (\ref{onedimbis}).

Thus, $z>0$ in $(0,\tau R)$. Furthermore, $z'(0)=0$, hence $z(0)>0$ from Hopf lemma. Another application of Hopf lemma shows that $z^{\prime}(\tau R)<0$. Therefore, there exists $\kappa>0$ such that $z>\kappa u_n$ in $[0,\tau R)$, whence
$$\frac{\gamma^*}{1+\kappa}\psi_1>u_n\ \hbox{ in }[0,\tau R),$$
which is a
contradiction with the definition of $\gamma^{\ast}$. 

Finally, we have obtained that $\mu<\lambda$, which means that $F_n(m,\tau)>F_1(2R,\tau)$.

We look for a reverse inequality. To that purpose, let $\varepsilon\in(0,1)$ and $R_0>0$ such that $\displaystyle
\frac{n-1}{R_0}<\varepsilon$. In the following computations, we always assume that $\tau R>R_0$. Define $u_n$ and $\lambda$ as
before. Let
$$\displaystyle \mu^{\prime}=\underline{\lambda}\left(\left(-\left(R-\frac{R_0}{\tau}\right),\left(R-\frac{R_0}{\tau}\right)\right), \tau(1-\varepsilon)\right)$$
and $w$ the normalized corresponding eigenfunction, so that
\[
\left\{
\begin{array}{l}
\displaystyle -w^{\prime\prime}(r)+\tau (1-\varepsilon) w^{\prime}(r)=\mu^{\prime}w(r) \displaystyle \mbox{ in }\left[0,R-\frac{R_0}{\tau}\right],\\
\\
\displaystyle w^{\prime}(0)=0,\ w>0\hbox{ in }\displaystyle\left[0,R-\frac{R_0}{\tau}\right),\ w\left(R-\frac{R_0}{\tau}\right)=0.
\end{array}
\right.
\] 
For all $R_0\leq x\leq \tau R$, define $\displaystyle v(x)=w\left(\frac{x-R_0}{\tau}\right)$, which satisfies
\[
\left\{
\begin{array}{l}
\displaystyle -v^{\prime\prime}(r)+(1-\varepsilon)v^{\prime}(r)=\frac{\mu^{\prime}}{\tau^2}v(r)\ \mbox{ in }\left[R_0,\tau
R\right],\\
\\
\displaystyle v^{\prime}(R_0)=0,\ v>0\hbox{ in }[R_0,\tau R),\ v(\tau R)=0.
\end{array}
\right.
\] 
Assume that $\lambda\geq \mu^{\prime}$. Since $(n-1)/R_0<\varepsilon$ and $u_n'(r)<0$ in $(0,\tau R]$, one therefore has
\[
\left\{
\begin{array}{ll}
\displaystyle -u_n^{\prime\prime}(r)+(1-\varepsilon)u_n^{\prime}(r)\geq \frac{\mu^{\prime}}{\tau^2}u_n(r) &\mbox{ in }[R_0,\tau
R],\\
\\
\displaystyle -v^{\prime\prime}(r)+(1-\varepsilon)v^{\prime}(r)=\frac{\mu^{\prime}}{\tau^2}v(r) &\mbox{ in }[R_0,\tau R].
\end{array}
\right.
\]
Arguing as before, we see that there exists $\gamma>0$ such that $\gamma u_n>v$ in $[R_0,\tau R)$. Define $\gamma^{\ast}$ ($>0$) as the
infimum of all such $\gamma$'s and define $z=\gamma^{\ast}u_n-v$, which is nonnegative in $[R_0,\tau R]$ and satisfies
$\displaystyle -z^{\prime\prime}+(1-\varepsilon)z^{\prime}-(\mu^{\prime}/\tau^2)z\ge 0$ in $[R_0,\tau R]$. 

Assume that $z(r)=0$ for some $r\in (R_0,\tau R)$. The strong maximum principle ensures that $z$ is $0$ in $[R_0,\tau R]$, which
means that $u_n=v$ in $[R_0,\tau R]$, which is impossible because $u_n^{\prime}(R_0)<0=v^{\prime}(R_0)$. 

Therefore, $z>0$ everywhere in $(R_0,\tau R)$. Furthermore, $z'(R_0)<0$, thus $z(R_0)>0$. On the other hand, by Hopf lemma, $z^{\prime}(\tau R)<0$. Thus, there exists $\kappa>0$ such that $z>\kappa v$ in $[R_0,\tau R)$, whence $(\gamma^*/(1+\kappa))u_n>v$ in $[R_0,\tau R)$. This contradicts the definition of $\gamma^{\ast}$.

Thus, we have established that $\lambda<\mu^{\prime}$. Straightforward computations (similar to those of the proof of (\ref{F1})) show
that
\[
\lambda<\mu^{\prime}=\frac{\tau^2}4\left(1-\varepsilon+\sqrt{(1-\varepsilon)^2-\frac{4\mu'}{\tau^2}}\right)^2
e^{-\sqrt{(1-\varepsilon)^2-\frac{4\mu'}{\tau^2}}(\tau R-R_0)},
\]
and, since $\lambda>F_1(2R,\tau)$, formula (\ref{F1}) and the fact that $m=\alpha_n R^n$ end the
proof of (\ref{Fntau}).\hfill\fin


\subsection{Proof of Proposition \ref{prosup}}\label{secminmax}

Let $m>0$ be fixed. For any $\Omega\in {\mathcal C}$ and any $v\in L^{\infty}(\Omega)$, one has the following min-max characterization of $\lambda_{1}(\Omega,v)$ (see \cite{bnv}):
\be\label{minmax}
\lambda_{1}(\Omega,v)=\sup_{\varphi}\ \inf_{\Omega} \left(\frac{-\Delta\varphi+v\cdot \nabla \varphi}{\varphi}\right),
\ee
where the supremum is taken over all functions $\varphi\in W^{2,n}_{loc}(\Omega)$ which are positive in $\Omega$. For $\varepsilon>0$, consider a domain $\Omega_{\varepsilon}\in {\mathcal C}$ included in $\left\{x=(x_{1},\ldots,x_{n})\in \R^n;\ \varepsilon<x_{1}<2\varepsilon\right\}$ and satisfying $\left\vert \Omega_{\varepsilon}\right\vert=m$, and let $v$ be any field in $L^{\infty}(\Omega_{\varepsilon},\R^n)$ with $\left\Vert v\right\Vert_{\infty}\leq \tau$. Define, for all $x=(x_{1},\ldots,x_{n})\in \Omega_{\varepsilon}$,
\[
\varphi_{\varepsilon}(x)=\sin\left(\frac{\pi x_{1}}{3\varepsilon}\right).
\]
Then, for all $x\in \Omega_{\varepsilon}$, 
\[
-\frac{\Delta\varphi_{\varepsilon}(x)}{\varphi_{\varepsilon}(x)}=\frac{\pi^2}{9\varepsilon^2},\quad \frac{v\cdot\nabla\varphi_{\varepsilon}(x)}{\varphi_{\varepsilon}(x)}\geq -\frac{2\tau\pi}{3\sqrt{3}\varepsilon}
\]
and thus
\[
\lambda_{1}(\Omega_{\varepsilon},v)\geq \frac{\pi^2}{9\varepsilon^2}-\frac{2\tau\pi}{3\sqrt{3}\varepsilon}.
\]
This shows that
\[
\sup_{\Omega\in{\mathcal{C}},\ \left\vert \Omega\right\vert=m} \underline{\lambda}(\Omega,\tau)\geq \frac{\pi^2}{9\varepsilon^2}-\frac{2\tau\pi}{3\sqrt{3}\varepsilon},
\]
and, since this is true for all $\varepsilon>0$, the assertion (\ref{supomega}) follows.\hfill\fin


\end{document}